  \documentclass{amsart}
  \usepackage[cp850]{inputenc}
  \RequirePackage{amssymb} % for the square for the proof environment.
  \RequirePackage{amsmath} % for multiple lines in equations
  \usepackage{amsmath}
  \usepackage[mathscr]{eucal}

  \usepackage{amsfonts}
  \usepackage{latexsym}
   \RequirePackage{amsbsy} % for Boldsymbol
   \RequirePackage{amsopn} % for DeclareMathOperator
   \RequirePackage{amsmath}% for multiple lines in equations
    \RequirePackage{amssymb}% for subsetneq

  \DeclareMathOperator{\calM}     {\mathcal M}%
  \DeclareMathOperator{\calW}     {\mathcal W}%
  \DeclareMathOperator{\calJ} {\mathcal J}%
  \DeclareMathOperator{\calU} {\mathcal U}%

   \newtheorem{thee}{Theorem}[section]
   \newtheorem{coor}[thee]{Corollary}
   \newtheorem{leem}[thee]{Lemma}
   
   \newtheorem{prro}[thee]{Proposition}
   
   \newtheorem{exxe}[thee]{Example}
   \newtheorem{reem}[thee]{Remark}

   \newcommand{\balf}
   {\renewcommand{\theenumi}{(\alph{enumi})}
   \renewcommand{\labelenumi}{\theenumi}
                        \begin{enumerate}}
  \newcommand{\ealf}   {\end{enumerate}
                        \renewcommand{\theenumi}{\arabic{enumi}}
                        \renewcommand{\labelenumi}{\theenumi.}}
  \newcommand{\bara}   {\renewcommand{\theenumi}{(\arabic{enumi})}
                        \renewcommand{\labelenumi}{\theenumi}
                        \begin{enumerate} }
  \newcommand{\eara}   {\end{enumerate}
                        \renewcommand{\theenumi}{\arabic{enumi}}
                        \renewcommand{\labelenumi}{\theenumi.}}

   \newcommand{\brom} {\renewcommand{\theenumi}{         (\roman{enumi})} \renewcommand{\labelenumi}{\theenumi}
                        \begin{enumerate} }
	            
  \newcommand{\erom}   {\end{enumerate}
                        \renewcommand{\theenumi}{\arabic{enumi}}
                        \renewcommand{\labelenumi}{\theenumi.}}

\begin{document}

   \title[ ]{NAGATA RINGS, KRONECKER
   FUNCTION RINGS AND RELATED SEMISTAR OPERATIONS}

  %AUTHOR INFO
  \author[ ]{Marco Fontana$^{1}$ \; and\ \; K. Alan Loper$^{2}$}
  % \address{Dipartimento di Matematica \\
  % Universit\`a degli Studi Roma Tre \\
  % Largo San Leonardo Murialdo, 1 \\ 
  % 00146 Roma, Italy}
  % \email{fontana@mat.uniroma3.it}
  % \author[ ]{ }
   %\address{Department of Mathematics
  % \\ Ohio State University-Newark\\ Newark, Ohio 43055 USA}
  % \email{lopera@math.ohio-state.edu}

  %\date{ } 
   
   \maketitle

  \dedicatory{\rm \begin{center} $^{1}$ Dipartimento di Matematica \\
  Universit\`a degli Studi Roma Tre \\
  Largo San Leonardo Murialdo, 1 \\ 
  00146 Roma, Italy\\ {\texttt{fontana@mat.uniroma3.it}} \\
  \texttt{              } \\ \rm $^{2}$ Department of Mathematics \\ 
  Ohio State
  University-Newark\\ Newark, Ohio 43055 USA \\
  {\texttt{lopera@math.ohio-state.edu}} \rm \end{center}}

  \bigskip
  \bigskip
  \bigskip

   \begin{abstract}

      In 1994, Matsuda  and Okabe introduced the notion of semistar 
  operation.
  This concept extends the classical concept of star operation (cf.
  for instance, Gilmer's book \cite{G}) and, hence, the related 
  classical
  theory of ideal systems based on the works by W.  Krull, E.  Noether, 
  H.
  Pr\"{u}fer and P.  Lorenzen from 1930's.

  In \cite{FL1} and \cite{FL2} the current authors investigated 
  properties 
  of the Kronecker function rings which arise from arbitrary semistar 
  operations on an integral domain $D$.  In this paper we extend that 
  study 
  and also generalize Kang's notion of a star Nagata ring 
  \cite{Kang:1987} 
  and \cite{Kang:1989} to the semistar setting.  Our principal focuses 
  are the 
  similarities between the ideal structure of the Nagata and Kronecker 
  semistar rings and between the natural semistar operations that these 
  two types of function rings give rise to on $D$.  
  \end{abstract}

  \bigskip

  \section{Introduction}

  A principal use of the classical star operations has been to construct
  Kronecker function
  rings   associated to an
   integral domain, in a more general context than the original one
   considered by L. Kronecker in 1882 \cite{K} (cf. \cite{Kr1},
   \cite{Kr2}, and  \cite{E} for a modern
   presentation of Kronecker's theory).  In this setting, one begins 
  with
   an integrally closed domain $\,D\,$ and a star operation $\,\star\,$ 
  on
  $\,D\,$  with the cancellation property known as e.a.b. (\sl endlich 
  arithmetisch
   brauchbar\rm ).  Then the Kronecker function ring is constructed as
  follows:

  \centerline{$
  {\textstyle\rm Kr}(D,\star) := \{f/g \; | \; \, f,g \in D[X] \setminus
  \{0\} \;\,
  \textrm{ and } \;\, \boldsymbol{c}(f)^\star \subseteq
  \boldsymbol{c}(g)^\star
  \} \, \cup \,\{0\}\,,
    $}
    
    \noindent (where $\,\boldsymbol{c}(h)\,$ denotes the content of a
    polynomial $\,h \in D[X]\,).$ \
  This domain turns out to be a B\'ezout overring of the polynomial ring
  $\,D[X]\,$
   such that $\,{\textstyle\rm Kr}(D,\star) \cap K = D\,$ (where $\,K\,$
  is the quotient field of $\,D\,$), cf.  \cite[Section 32]{G}.

  In 1994, Okabe and Matsuda \cite{OM1} introduced the more ``flexible''
  notion of
  semi\-star operation $\,\star\,$ of an integral domain $\,D\,,$ \ as a
  natural generalization of the notion of star operation, allowing $\, D
  \neq D^\star\,$ (the definition is given in Section 2); cf. also
  \cite{MS}, \cite{MSu} and \cite{OM2}.  In several recent
  papers, the classical construction introduced by Kronecker has been
  further
  generalized so that we can begin with {\sl any} integral domain (not
  ne\-cessa\-rily integrally closed) $\,D\,$ and {\sl any} semistar 
  operation
  (not necessarily e.a.b.) $\,\star\,$ of $\,D\,$ and, in a natural
  manner, construct a Kronecker function ring, still denoted here by
  $\,{\textstyle\rm Kr}(D,\star)\,,$ \ which preserves the main 
  properties
  of the ``classical'' Kronecker function ring (axiomatized by 
  Halter-Koch
  \cite{HK3}, see also \cite{ADF:1987}), cf.  the works by Okabe and
  Matsuda \cite{OM3}, by Matsuda \cite{M}, and by Fontana and Loper
  \cite{FL1}, \cite{FL2}.  \ One evidence of the ``naturalness'' of this
  construction is that this general Kronecker function ring gives rise 
  to
  an e.a.b. semistar operation $\,\star_a\,$ (which can be 
  ``restricted''
  to an e.a.b. star operation, denoted by $\,\dot{\star}_a\,,$ \ of the
  integrally closed overring $\,D^{\star_{a}}\,$ of $\,D\,$) and then
  $\,{\textstyle\rm Kr}(D,\star)\,$ can be viewed in the classical star
  e.a.b. setting by use of $\,\star_a\,$ (or, more precisely, by use of
  $\,\dot{\star}_a\,)$.

  Another overring of $D[X]$ which has been much studied is the Nagata
  ring of $D$, i.e.  $D(X) := \{f/g \; | \; \, f,g \in D[X] \; \; 
  \mbox{\rm and}
  \; \; {\boldsymbol c}(g) = D\}\,, $ \ cf.  \cite{Krull}, \cite[page 
  27]{Samuel}, \cite[page
  18]{Nagata}, \cite[Section
  33]{G}, \cite[Chapter IV]{Huckaba}.  The interest in $\,D(X)\,$ is due
  to the fact that this ring has some ``nice'' pro\-per\-ties that 
  $\,D\,$
  itself need not have, mantaining in any case a strict relation with 
  the
  ideal structure of $\,D\,,$ \ (for instance, for each ideal $\,I\,$ of
  $\,D\,$, we have $\, ID(X) \cap D = I\,$ and $\, D(X)/ID(X) \cong
  (D/ID)(X)\,$).  \ Among the ``new'' properties acquired by $\,D(X)\,$ 
  we
  mention\ (a) the residue field at each maximal ideal of $\,D(X)\,$ is
  infinite;\ (b) an ideal contained in a finite union of ideals is
  contained in one of them \cite {QB};\ (c) each finitely generated
  locally principal ideal is principal \cite{A:1977}.  Furthermore, the
  canonical map \ Spec$(D(X)) \rightarrow $ Spec$(D)\,$ is a 
  homeomorphism
  if and only if the integral closure of $\,D\,$ is a Pr\"ufer domain
  \cite{ADF}.  The relation between the Nagata ring and the Kronecker
  function ring was investigated by Arnold \cite{Arnold:1969}, by Gilmer
  \cite{Gilmer:1970} and by Arnold and Brewer \cite{Arnold-Brewer:1971}.
  In particular, if $\,b\,$ denotes the star operation of $\,D\,$ 
  defined
  on the fractionary ideals $\,I\,$ of $\,D\,$ by $ I \mapsto I^b := 
  \cap
  \{ IV \;| \;\, V \; \mbox{is a valuation overring of} \; D \}\,,$ \ 
  then
  Arnold \cite[Theorem 7]{Arnold:1969} proved that D is a Pr\"ufer 
  domain
  if and only if $\,D(X) = {\textstyle\rm Kr}(D, b)\,.$

   A generalization of the Nagata ring construction was
  considered by Kang \cite{Kang:1987}, \cite{Kang:1989}: for each star
  operation $\,\star\,$
  of $\,D\,$ he studied the ring $\,\{f/g \; | \; \, f, g \in D[X] \;
  \,\mbox{and} \; \,$ $ {\boldsymbol c}(g)^\star =D\}\,.$ \ In 
  particular,
  Kang
  proved that, \sl mutatis mutandis\rm, many properties of the
  ``classical''
  Nagata ring still hold in this more general context.

  \vskip 4pt

  In the present paper we further generalize the previous construction 
  so
  that,
  given any domain $\,D\,$ and any semistar operation $\,\star\,$ on
  $\,D\,,$ \ we define the semistar Nagata ring as follows:

  \centerline{$
  {\textstyle\rm Na}(D,\star) :=\{f/g \; | \;\, f, g \in D[X] \;\;
  \mbox{and} \;\; {\boldsymbol c}(g)^\star =D^\star\}\,.
  $}

  \noindent We then study the ideal structure of $\,{\textstyle\rm
  Na}(D,\star)\,$ and compare it to that of $\,{\textstyle\rm
  Kr}(D,\star)\,.$ \ We also show how $\,{\textstyle\rm Na}(D,\star)\,$
  gives rise to a very natural semi\-star operation, denoted by
  $\,\tilde{\star}\,,$ \ which plays a role analogous to that of the
  semistar operation $\,\star_a\,$ in the Kronecker setting.  In the
  star operation case, $\,\tilde{\star}\,,$ coincides with the operation
  $\,\star_{w}\,$ considered recently by D.D. Anderson and S.J. Cook
  \cite{Anderson/Cook:2000}.

  In Section 2 we give some background information concerning semistar
  operations and some preliminary results concerning the class of
  quasi--$\star$--ideals (i.e. ideals $\,I\,$ such that $\,I^\star \cap 
  D
  = I\,$), a
  more general class than that of $\,\star$--ideals (i.e. ideals $\,I\,$
  such that $\,I^\star  = I\,$),
  which plays an
  important role when $\,\star\,$ is a semistar ope\-ra\-tion.

  In the third section we define and study the
   semistar  Nagata rings.  For instance, we show that there is a 
  natural
  1-1
  correspondence between the maximal ideals of $\,\mbox{\rm
  Na}(D,\star)\,$ and the maximal elements in the set of all proper
  quasi--$\star$--ideals of $\, D\,;$ \ in particular, the 
  $\,t$--maximal
  ideals of an integral domain $\,D\,$ are all obtained as contractions 
  to
  $\,D\,$ of the maximal ideals of $\,\mbox{\rm Na}(D, v)\,.$ \ 
   We prove also that $\,\mbox{\rm Na}(D, \star)
  =\mbox{\rm Na}(D,\tilde{\star})\,$.  \ Furthermore, we
  show
  that there is a strict link between the semistar operation
  $\,\tilde{\star}\,,$ \
   the maximal elements $\,P\,$ in the set of all proper
  quasi--$\star$--ideals of $\, D\,$ and the valuation overrings of
  $\,D_{P}\,.$ \ More precisely, if we say that a $\,\star$--valuation
  overring of $\,D\,$ is a valuation overring of $\,D\,$ such that 
  $\,F^\star
  \subseteq FV\,$ for each finitely generated fractionary ideal $\,F\,$
  of $\,D\,,$ \ then we show that a valuation overring $\,V\,$ of 
  $\,D\,$
  is a $\,\tilde{\star}$--valuation overring of $\,D\,$ if and only if
  $\,V\,$ is an overring of $\,D_P\,,$ \ for some $\,P\,$ maximal in the
  set
  of all proper quasi--$\star$--ideals of $\, D\,.$

  In the fourth section we
  recall from \cite{FL1} some results concerning $\,{\textstyle\rm
  Kr}(D,\star)\,$ and $\,\star_a\,$ and examine the interplay with
  $\,{\textstyle\rm Na}(D,\star)\,$ and $\,\tilde{\star}\,.$  \ In
  particular, we
  show that each maximal element $\,Q\,$ in the set of all proper
  quasi--$\star_{a}$--ideals of $\, D\,$ is determined uniquely by a
  $\,\star$--valuation overring of $\,D\,$ (dominating $\,D_{Q}\,$) and
  there is a natural 1-1 correspondence between the maximal ideals of
  $\,\mbox{\rm Kr}(D,\star)\,$ and the minimal $\,\star$--valuation
  overrings of $\,D\,$.

  In the final section we examine more closely the relationship between
  $\,\tilde{\star}\,$ and $\,\star_a\,$ and we show that it is hopeless 
  to
  try to attain an equality by applying
  $\,\widetilde{(\mbox{-})}\,$ and $\,(\mbox{-})_a\,$ in different 
  orders
  to an arbitrary semistar operation.

  \vskip 6pt

  We use Gilmer's book \cite{G} as our main reference.  Any unexplained
  material is as in \cite{G} and \cite{Ka}.  Many preliminary
  results on semistar operations and applications  appear in
  conference proceedings   (in
  particular, \cite{FL1} and \cite{FL2}), and hence are not easily 
  available.  Because of this possible hindrance,  we  briefly
  restate the 
  principal
  definitions and statements of the main properties that we will 
  need so that
   the present work will be self--contained.

  \noindent Note that the ``module systems'' approach, developed very 
  recently by
  Halter-Koch in \cite{Koch:JA}, provides a general setting for
  (re)considering semistar operations and, in particular, many of the
  constructions related to the semistar operations considered in the
  present paper.  However, since many background results of our paper 
  are
  proved in an earlier work by Fontana and Huckaba \cite{FH}, which has
  inspired and provided the foundation also of \cite{Koch:JA}, we 
  maintain
  the level of generality of this paper within the more classical 
  ``semistar''
  setting.

  \vskip 4pt

  \section{Background and preliminary results}

  For the duration of this paper $\,D\,$ will represent an integral 
  domain
  with quotient field $\,K\,$.  Let $\,\boldsymbol{\overline{F}}(D)\,$
  represent the set of all nonzero $\,D$--submodules of $\,K\,$.  Let
  $\,\boldsymbol{F}(D)\,$ represent the nonzero fractionary ideals of
  $\,D\,$ (i.e. $\,E \in \boldsymbol{\overline{F}}(D)\,$ such that $\,dE
  \subseteq D\,$ for some nonzero element $\,d \in D\,$).  Finally, let
  $\,\boldsymbol{f}(D)\,$ represent the finitely generated
  $\,D$-submodules of $\,K\,$.

  A mapping $\,\star : \boldsymbol{\overline{F}}(D) \rightarrow
  \boldsymbol{\overline{F}}(D)\,$, $\,E \mapsto E^\star\,$
  is called \it a semistar operation of $\,D\,$ \rm if, for all $\,z \in
  K\,$,
  $\,z \not =
  0\,$ and for all $\,E,F \in \boldsymbol{\overline{F}}(D)\,$, the
  following
  properties hold:

  $\mathbf{ \bf (\star_1)} \; \; (zE)^\star = zE^\star \,; $

  $\mathbf{ \bf (\star_2)} \; \;   E \subseteq F \;\Rightarrow\; E^\star
  \subseteq
  F^\star \,;$

  $\mathbf{ \bf (\star_3)} \; \;   E \subseteq E^\star \; \; \textrm { 
  and
  }
  \; \; E^{\star \star} := (E^\star)^\star = E^\star\,. $

  %REMARK 2.1
  \begin{reem} \label{rk:2.1} \rm  Let $\,\star\,$ be a semistar
  o\-pe\-ra\-tion of $\,D\,$.

  \bf (a) \rm  If $\,\star\,$ is a semistar operation such that 
  $\,D^\star
  =
  D\,$, then the map $\,\star :
  \boldsymbol{F}(D) \rightarrow \boldsymbol{F}(D)\,$, \, $E \mapsto
  E^\star\,,$ is called a \it
  star operation of \rm $\,D\,.$ \ Recall from \cite[(32.1)]{G} that a
  star
  operation $\,\star\,$ verifies the properties $\,\mathbf{ \bf
  (\star_2)}\,, \mathbf{ \bf (\star_3)}\,,$ for all $\,E, F \in
  \boldsymbol{F}(D)\,;$ moreover, the property $\,\mathbf{\bf
  (\star_1)}\,$ can be restated as follows: for each $\,z\in K\,, \ z 
  \neq
  0\,$ and for each $\,E \in \boldsymbol{F}(D)\,,$

   $\mathbf{(\star\star_1)}\; \; $ $(zD)^\star = zD\,, \; \; 
  (zE)^{\star}
  = zE^{\star}\,.$

  \noindent If $\,\star\,$ is a semistar operation of $\,D\,$ such that
  $\,D^\star = D\,,$ then we will write
  often in the sequel that $\,\star\,$ is \it a (semi)star
  operation of $\,D\,,$ \rm to
  emphasize the fact that the semistar operation $\,\star\,$ is an
  extension to
  $\,\boldsymbol{\overline{F}}(D)\,$ of \it a ``classical'' star 
  operation
  \rm $\,\star\,,$ i.e. a map
  $\, \star : \boldsymbol{F}(D)\rightarrow \boldsymbol{F}(D)\,,$ 
  verifying
  the properties
  $\,\mathbf{(\star\star_1)}\,, \mathbf{(\star_2)}\,$ and $\,
  \mathbf{(\star_3)}\,$ \cite[Section 32]{G}.  Note that
  not every semistar operation is an extension of a star operation
  \cite[Remark 1.5 (b)]{FH}.

  \bf (b) \rm  The \it trivial semistar operation on $\,D\,$ \rm is the
  semistar
  operation constant onto $\,K\,$, i.e. the semistar operation
  $\,\star\,$ such that $\, E^\star = K\,,$ for each  $\, E \in
  \boldsymbol{\overline{F}}(D)\,$. \ Note that  $\,\star\,$ is the
  trivial semistar operation on $\,D\,$ if and only if $\, D^{\star} =
  K\,$. (As a matter of fact, if $\,D^{\star}=K\,$,  then for each $\, E
  \in
  \boldsymbol{\overline{F}}(D)\,$ and for each $\,e \in E\,,\ e \neq
  0\,,$\ we have  $\,eD \subseteq E \,$ and thus $\, K =eK =eD^{\star}
  \subseteq
  E^{\star} \subseteq K\,.$)

   \bf (c) \rm  Let $\,D\,$ be an integral domain and $\,T\,$ an
  overring
  of $\,D\,$.  Let $\,\star\,$ be a semistar operation of $\,D\,$ and
  define
    $\dot{\star}^{\mbox{\tiny \it \tiny T}}: 
  \overline{\boldsymbol{F}}(T)
    \rightarrow \overline{\boldsymbol{F}}(T)$ by setting:

  \centerline{$
  E^{\dot{\star}^{\mbox{\tiny \it \tiny T}}} := E^\star\,, \;
  \mbox{ for each } \;  E \in  \overline{\boldsymbol{F}}(T) \ 
  (\subseteq \overline{\boldsymbol{F}}(D))\,.
  $}

  \noindent Then, we know \cite[Proposition 2.8]{FL1}:

  \hskip 0.2cm \bf (c.1)  \rm The operation $\,\dot{\star}^{\mbox{\tiny
  \it
  T}}\,$ is a semistar operation of $\,T\,.$

  \hskip 0.2cm \bf (c.2) \rm When $\,T = D^\star $,\ then
  ${\dot{\star}}^{{\mbox{\tiny \it D}}^{\star}}$ defines a (semi)star
  operation of $D^\star\,.$

   \bf (d) \rm If $\,\star_{1}\,$ and $\,\star_{2}\,$ are two semistar
  operation
  of $\,D\,$,\ we say that $\; \star_{1} \leq \star_{2}\;$ if
  $\,E^{\star_{1}} \subseteq E^{\star_{2}}\,$, \ for each $\,E \in
  \boldsymbol{\overline{F}}(D)\,;$ \ in this case,
  $\,(E^{\star_{1}})^{\star_{2}} = E^{\star_{2}}\,.$
  \end{reem}

  \vspace{-2pt}

  We refer to the collection
  % 
  %\vspace{.1in}
  % \hskip 2.6cm
  $ \boldsymbol{\overline{F}}^\star (D) := \{E^\star \; | \;\, E
  \in
  \boldsymbol{\overline{F}}(D) \} $
  %\vspace{.1in}
  % 
  % \noindent
  [respectively,  
  %\hskip 0.92 cm
  $\boldsymbol{F}^\star (D) := \{I
  \in \boldsymbol{F}(D)\; |\;\, I = H^\star \,\textrm { with } \, H \in
  \boldsymbol{F}(D) \}\,;$
  % 
  %\vspace{.1in}
  % \hskip 2.6cm
  $ \boldsymbol{f}^\star (D) := \{J \in \boldsymbol{F}(D)\;
  |\;\, J = F^\star \, \textrm { with } \, F \in \boldsymbol{f}(D) \} 
  \,$] \
  % 
  % %\vspace{.1in}
  % 
  \noindent as
  %\begin{itemize}
  %\item
  \it the $\,\star$--$D$--submodules of $\,K$\, \rm
  %\item
  [respectively, \, \it the (fractionary) $\,\star$--ideals of $D\ ;\;$
  %\item
  the (fractionary) $\,\star$--ideals of $\,D\,$ of finite type\rm \,].

  %\end{itemize}

  These labels seem natural, but can be problematic.  As a matter of 
  fact,
  if $\,I \in \boldsymbol{F}(D)\,,$\ then $\,I^\star\,$ is not
  necessarily a fractionary ideal of $\,D\,$ and so it does not 
  necessarily belong to $\,\boldsymbol{F}^\star (D)\,$ (e.g. if $\,(D: 
  D^\star)=
  0\,,$\ then $\,D^\star \not\in\boldsymbol{F}^\star (D)\,$). 
 For instance
  if $\,T\,$
  is an overring of an integral domain $\,D\,$ such that the conductor
  $\,(D:T) =0\,$ and if $\,\star := {\star}_{\mbox{\tiny\{\it \tiny
  T\}}}\,$ is the semistar operation of $\,D\,$ defined by 
  $\, E^{\star_{\mbox{\tiny \{\it
  \tiny
  T\}}}}:= ET\,$ \ for
  each $\,E\,
  \in \overline{\boldsymbol{F}}(D),$
   then it is easy to see that $\,\boldsymbol{F}^{\star} (D)\,$ is 
  empty.
  %will consist of only the zero ideal.
  So we need a more general notion than $\,\star$--ideal, when
  $\,\star\,$ is a semistar
  ope\-ra\-tion.

  \vspace{4pt}

  %\begin{deef}
   \rm  Let $\,I \subseteq D\,$ be a nonzero ideal of $\,D\,$ and let
  $\,\star\,$ be a semistar operation on $\,D\,$.  We say that $\,I\,$ 
  is
  \it a
  quasi--$\star$--ideal of $\,D\,$ \rm if $\,I^\star \cap D = I\,$. \
  Similarly, we
  designate by \it  quasi--$\star$--prime \rm [respectively,\, \it 
  $\,\star$--prime
  \rm] of $D$ a quasi--$\star$--ideal [respectively,\, an integral
  $\,\star$--ideal]  of $\,D\,$  which is also a prime ideal. We 
  designate by \it
   quasi--$\star$--maximal \rm
  [respectively,\, \it  $\,\star$--maximal \rm]  of $\,D\,$  a maximal
  element in the
  set of all proper quasi--$\star$--ideals [respectively,\, integral
  $\star$--ideals] of $\,D\,.$
  %\end{deef}

  \vspace{4pt}

  Note that if $\,I \subseteq D\,$ is a $\,\star$--ideal, it is also a
  quasi--$\star$--ideal and, when $\,D = D^\star\,$
  the notions of quasi--$\star$--ideal and
  integral $\,\star$--ideal coincide.  When $\,D \subsetneq D^\star
  \subsetneq K\,$, we can ``restrict'' the semistar operation 
  $\,\star\,$
  on
  $\,D\,$ to the nontrivial (semi)star operation on $\,{D^\star}\,$,
  denoted by $\,{\dot{\star}}^{{\mbox{\tiny \it D}}^{\star}}\,$, \ or
  simply by $\,{\dot{\star}}\,$, and defined in Remark~\ref{rk:2.1} (c),
  and we have a strict relation between the quasi--$\star$--ideals of
  $\,D\,$ and
  the $\,{\dot{\star}}$--ideals of $\,{D^\star}\,$.
  \vskip -16pt
  %LEMMA 2.2
  \begin{leem}\label{lm:1}
   \sl Assume the notation of the preceding paragraph.
  Then:
   
  \centerline{$ I
  \, \mbox{ is a quasi--$\star$--ideal of $D$}\, \Leftrightarrow \, I = 
  L
  \cap D\,, \; \mbox{where} \; L \subseteq D^{\star} \, \mbox{ is
  a $\dot{\star}$--ideal of $D^\star$}\,.$}
  \end{leem} \rm
  \vskip -3pt  \noindent {\bf Proof}.  The proof follows  from $
   I^\star \cap D \subseteq (I^\star \cap D)^\star \cap D \subseteq
  I^{\star \star} \cap D^\star \cap D = I^\star \cap D$.
  \hfill $\Box$

  \vspace{4pt}

  Note that this also gives a means of constructing 
  quasi--$\star$--ideals
  and, in particular, quasi--$\star$--ideals containing a given ideal.  
  If
  $I
  \subseteq D$ is a nonzero ideal of $D$, then $\,I^\star \cap D\,$ is a
  quasi--$\star$--ideal of $D$ which contains $I$.

  \noindent We denote by\ Spec$^\star(D)$\ [respectively, 
  Max$^\star(D)$;
  QSpec$^\star(D)$; QMax$^\star(D)$] the set of all $\star$--primes
  [respectively,\, $\star$--maximals;\ quasi--$\star$--primes;\
  quasi--$\star$--maximals] of $\,D\,$.

  \vspace{4pt}

  As in the classical star-operation setting, we associate to a semistar
  ope\-ra\-tion $\,\star\,$ of $\,D\,$ a new semistar operation
  $\,\star_f\,$ as follows.
  %\begin{deef}
   \rm Let $\,\star\,$ be a semistar operation of a domain $\,D\,$.  \ 
  If
   $\,E \in \boldsymbol{\overline{F}}(D)\,$ we set
  % 
  % \centerline{
  $
  E^{\star_f} := \cup \{F^\star \;|\;\, F \subseteq E,\, F \in
  \boldsymbol{f}(D)
  \}\,.
  $
  % }

  \noindent We call $\,\star_f\,$ \it the semistar operation of finite
  type of
  $D$ \rm associated to $\,\star\,.$ \ If $\,\star = \star_f\,$,\ we say
  that
  $\,\star\,$ is \it a semistar operation of finite type of $\,D\,$.  
  \rm
  \ Note that $\,\star_f \leq \star\,$ and $\,(\star_f)_f = \star_f\,$,\
  so
  $\,\star_f\,$ is a semistar operation of finite type of $\,D\,.$
  \, For instance, if $\,v\,$ is \it the  $v$--(semi)star operation on 
  $D$
  \rm defined by $E^v := (E^{-1})^{-1},$ for each $E \in
  \overline{\boldsymbol{F}}(D)\,,$ with $E^{-1} := (D :_{\mbox{\tiny \it
  K}} E) := \{ z \in K \; | \;\; zE \subseteq D \}\,$ \cite[Example
  1.3 (c) and Proposition 1.6 (5)]{FH},\, then the semistar operation of
  finite type $\,v_{f}\,$ associated to $\,v\,$ is called \it the
  $t$--(semi)star
  operation on $D$\, \rm (in this case $D^v =D^t = D$).

  Both the Kronecker function rings and the Nagata rings considered in 
  the present paper are defined in a natural way for a general semistar 
  operation.  A principal theme of the paper is that both of 
  these classes of rings can be recast as Kronceker function rings and 
  Nagata rings of certain natural semistar operations of finite type.  
  So the entire 
  theory could be stated in terms of semistar operations of finite 
  type.  It seems worthwhile to us to keep the more general setting so 
  that, for example, we can talk about the Kronecker function ring and 
  the Nagata ring associated to the classical $v$ operation (which is 
  rarely of finite type).

  %LEMMA 2.3

  \begin{leem}\label{lm:2} \sl Let $\,\star\,$ be a semistar operation 
  of
  an integral
  domain $\,D\,$. \ Assume that $\,\star\,$ is not trivial and that
  $\,\star =
  \star_f\,$. \ Then

  %\begin{enumerate}
  \bara
  \rm \bf \item \sl Each proper quasi--$\star$--ideal is contained in a
   quasi-$\star$-maximal.
  \rm \bf \item \sl Each quasi--$\star$--maximal is a
  quasi--$\star$-prime.

  \rm \bf \item \sl If
   $\, Q\, \mbox{ is a quasi--$\star$--maximal ideal of $D$}\,$ then
   $\, Q = M \cap D\,, \; $ for some $\,\dot{\star}$--maximal ideal $\,
  M\,$
   of $D^\star\,.$
   
   \rm \bf \item \sl If $\, L \subseteq D^{\star} \,$ is a 
  $\dot{\star}$--prime ideal of
  $\,D^\star,$ \ then $\,L\cap D\,$ is a quasi--$\star$--prime ideal of
  $\,D\,.$
   \rm \bf \item \sl Set
  \;\; $
  \Pi^\star : = \{P \in \mbox{\rm Spec}(D)\; |\;\, P \not = 0 \, \mbox{
  and } \, P^\star \cap D \not = D\}\,.
  $

  \noindent Then \,\rm QSpec$^\star(D) \subseteq \Pi^\star$\, \sl and 
  the
  set of maximal elements
  of $\;\Pi^\star,$\ denoted by $\; {\Pi^\star_{\mbox{\rm \tiny
  max}}} ,$ \ is nonempty and coincides with \rm QMax$^\star(D)$.

  \eara
  %\end{enumerate}

  \end{leem} \rm

  \vskip -3pt \noindent {\bf Proof}. The proof is straightforward.
 \hfill  $\Box$

  \vspace{5pt}

  Note that, in general, the restriction to $\,D\,$ of a
  $\dot{\star}$--maximal ideal of $D^\star$  is a
  quasi--$\star$--prime ideal of $\,D\,,$ \ but not necessarily a
  quasi--$\star$--maximal ideal of $\,D\,,$ \ and if $\, L\,$ is an 
  ideal
  of
  $\,D^\star\,$ and $\,L\cap D\,$ is a quasi--$\star$--prime ideal of
  $\,D\,,$ then $\, L\,$ is not necessarily a $\dot{\star}$--prime ideal
  of
  $\,D^\star,$ (cf.  the Remark~\ref{new}).

  For the sake of simplicity, when $\,\star = \star_f\,$, \ we will 
  denote
  simply by \, ${\calM}(\star) \,$, the nonempty set
  $\,\Pi^\star_{\mbox{\rm \tiny max}} = \mbox{QMax}^\star(D)\, $.

  \vspace{5pt}

  If $\,\Delta\,$ is a nonempty set of prime ideals of an integral 
  domain
  $\,D\,$, \ then the semistar operation $\,\star_\Delta\,$ defined
  on $\,D\,$ as follows

  \centerline{$
  E^{\star_\Delta} := \cap \{ED_P \;|\;\, P \in \Delta\}\,,
  \; \; \textrm {  for each}    \; E \in \boldsymbol{\overline{F}}(D)\,,
  $}

  \noindent is called \it the spectral semistar operation associated to
  \rm
  $\,\Delta\,$. \ If $\,\Delta = \emptyset\,$, \ then we can extend the
  previous
  defintion by setting $\,E^{\star_\emptyset} := K\,$, \ for each $\,E 
  \in
  \boldsymbol{\overline{F}}(D)\,$, \ i.e. $\,\star_\emptyset\,$ is the
  trivial
  semistar operation on $\,D\,$ (constant onto $\,K\,$; \ cf. 
  Remark~\ref{rk:2.1} (b)).

  %LEMMA 2.4
  \begin{leem}\label{lm:3} \sl Let $\,D\,$ be an integral domain and let
  $\; \emptyset \not = \Delta
  \subseteq \mbox{\rm Spec}(D)\,$. \ Then:

  %\begin{enumerate}
  \bara

  \rm \bf \item \sl  $E^{\star_\Delta}D_P = ED_P$, \, for each $E \in
  \boldsymbol{\overline{F}}(D)\,$ and for each $\,P \in \Delta\,$.
  \rm \bf \item \sl  $(E \cap F)^{\star_\Delta}= E^{\star_\Delta} \cap
  F^{\star_\Delta}$,\, for all $E,F \in \boldsymbol{\overline{F}}(D)$.
  \rm \bf \item \sl  $P^{\star_\Delta} \cap D = P$, \, for each $P \in
  \Delta$.
  \rm \bf \item \sl If $\,I\,$ is a nonzero integral ideal of $\,D\,$ 
  and
  $\,I^{\star_\Delta} \cap D \not = D\,$\, then there exists $\,P \in
  \Delta\,$
  such that $\,I \subseteq P\,$.

   \rm \bf \item \sl Assume that the set of
  maximal elements
  $\,\Delta_{\mbox{\rm\tiny max}}\,$ of $\,\Delta\,$ is also nonempty 
  and
  that each $\,P \in \Delta\,$ is contained in some $\,Q \in
  \Delta_{\mbox{\rm \tiny max}}\,$.  Then:

  \centerline{$\star_\Delta = \star_{\Delta_{\mbox{\rm\tiny max}}}\,. $}

  %\end{enumerate}
  \eara

  \end{leem} \rm

  \vskip -3pt \noindent {\bf Proof}: \cite[Lemma 4.1]{FH} and, for 
  (5),
  \cite[Remark 4.5]{FH}.  \hfill $\Box$

  \vspace{5pt}

  A semistar operation $\,\star\,$ of an integral domain $\,D\,$ is 
  called
  \it a
  spectral semistar operation \rm if there exists $\,
  \Delta
  \subseteq \mbox{\rm Spec}(D)\,$ such that $\,\star =
  \star_\Delta\,$. \ We say that
  $\,\star\,$ \it posesses enough primes \rm or that $\,\star\,$ is \it 
  a
  quasi-spectral
  semistar operation of $\,D\,$ \rm if, for each nonzero ideal $\,I\,$ 
  of
  $\,D\,$ such that
  $\,I^\star \cap D \not = D\,$, \ there exists a quasi--$\star$--prime
  $\,P\,$ of
  $\,D\,$ such that $\,I \subseteq P\,$. \ Finally, we say that
  $\,\star\,$ is \it
  a stable semistar operation on $\,D\,$ \rm if \
  $
  (E \cap F)^\star = E^\star \cap F^\star, \, \textrm {for all} \,
  E,F \in \boldsymbol{\overline{F}}(D) \,.
  $

  \vspace{4pt}

  %LEMMA 2.5
  \begin{leem}\label{lm:4}
   \sl Let $\,\star\,$ be a nontrivial semistar operation of
  an integral domain $D$.
  \bara
  \bf \item \sl   $\,\star\,$ is spectral if and only if
  $\,\star\,$ is quasi-spectral and stable.
  \bf \item \sl Assume that $\,\star = \star_f\,$. \ Then
  $\,\star\,$ is quasi-spectral and $\,\calM(\star) \not = 
  \emptyset\,$. 
  \eara\end{leem} \rm

  \vskip -3pt \noindent {\bf Proof}.  (1)  The ``only if" part is a
  consequence of Lemma~\ref{lm:3} (2) and (4).  The ``if" part is proved
  in
  \cite[Theorem 4.12 (3)]{FH}.\  (2)   is a restatement of 
  Lemma~\ref{lm:2}. \hfill $\Box$

  \vspace{4pt}

  If $\,\star\,$ is a semistar operation of an integral domain $D$ and 
  if
  $\,\Pi^\star \not = \emptyset\,$, the nontrivial semistar operation

  \vskip -3pt \centerline{$ \star_{sp} := \star_{\Pi^\star} $}

  \noindent is called \it the spectral semistar operation associated to
  $\,\star\,$.
  \rm

  %LEMMA 2.7
  \begin{leem} \label{lm:7} \sl Let $\,\star\,$ be a nontrivial semistar
  operation.
  \bara
  \bf \item \sl  $\,\star\,$ is spectral if and only if $\,\star = 
  \star_{sp}\,$.

   \bf \item \sl  Assume that $\,\Pi^\star \not = \emptyset\,$.  Then
  the
  following statements are equivalent:
  \eara
  %\begin{enumerate}
 {\brom
 \rm \bf \item \sl $\star_{sp} \leq \star\,$;

   \rm \bf \item \sl $\star$ is quasi-spectral;

  \rm \bf \item \sl $E^\star = \cap \{E^\star D_P \;|\;\, P \in
  \Pi^\star\}\,$,
  $\;$ for each $\,E \in \boldsymbol{\overline{F}}(D)\,$.

  %\end{enumerate}
  \erom }

  \end{leem}\rm

  \vskip -3pt \noindent {\bf Proof}: (1) 
   \cite[Corollary 4.10]{FH};\,  (2) 
   \cite[Proposition 4.8]{FH}. \hfill $\Box$

  %COROLLARY 2.8
  \begin{coor}\label{cor:9} \sl Let $\,\star\,$ be a nontrivial semistar
  operation.  Set

   \centerline{$ \tilde{\star} := (\star_f)_{sp} \,.$}

  %\noindent Then:

  %begin{enumerate}
  \bara

  \rm \bf \item \sl $\tilde{\star} = \star_{{\calM}(\star_f)} \leq
  \star_f\,$\
  (in particular, $\,\tilde{\star}\,$ is not trivial).
  \rm \bf \item \sl For each $\,E \in \boldsymbol{\overline{F}}(D)\,$,

  %\begin{itemize}
  \balf
  \rm \bf \item \sl $E^{\star_f} = \cap \{E^{\star_f}D_Q \;|\;\, Q \in
  {\calM}(\star_f) \}\,$;

  \rm \bf \item \sl $E^{\tilde{\star}} = \cap \{ED_Q \;|\;\, Q \in
  {\calM}(\star_f) \}\,$.

  \ealf

  %\end{itemize}

  %\end{enumerate}
  \eara

  \end{coor} \rm

  \vskip -3pt \noindent {\bf Proof}. (1) is a consequence of 
  Lemma~\ref{lm:4} (1), 
  Lemma~\ref{lm:3} (5), Lemma~\ref{lm:7} (2) and Lemma ~\ref{lm:2} (3).
  (2).  The first equality follows from Lemma ~\ref {lm:7} (2) ((ii)
  $\Rightarrow$
  (iii)), Lemma~\ref{lm:4} (1) and Lemma~\ref{lm:2} (3).  The second
  equality follows from
  (1) and from the definition of spectral semistar operation. \hfill
  $\Box$

  %\vspace{4pt}

  %REMARK 2.9
  \begin{reem} \label{rem:2.12} \rm {\bf (a)} Note that, when 
  $\,\star\,$
  is the (semi)star
   $v$--operation, then the (semi)\-star operation $\,\tilde{v}\,$
  coincides with
   \it the (semi)star operation $\,w\,$ \rm defined as follows:
    
  \vskip -3pt  \centerline{$E^{w} :=
  \cup \{(E:H)\; | \;\,H \in \boldsymbol{f}(D) \mbox{ and } H^v = D 
  \}\,,
   \;\; \mbox{ for each } E \in \boldsymbol{\overline{F}}(D)\,.$}
   
       \noindent This (semi)star operation was first considered by J. 
  Hedstrom and
  E. Houston in 1980 \cite[Section 3]{Hedstrom/Houston: 1980} under the
  name of F$_{\infty}$--operation.  Later, starting in  1997, this 
  operation was
       intensively studied by W. Fanggui and R. McCasland (cf.
       \cite{Fanggui:1}, \cite{Fanggui:2}, 
  \cite{Fanggui/McCasland:1997},
       \cite{Fanggui/McCasland:1999}, ) under the name of 
  $w$--operation.
  Note
       also that the notion of $w$--ideal coincides with the notion of
       semi-divisorial ideal considered by S. Glaz and W. Vasconcelos in
  1977
       \cite{Glaz/Vasconcelos:1977}.  Finally, in 2000, for each
  (semi)star
       operation $\,\star\,$, D.D. Anderson and S.J. Cook
       \cite{Anderson/Cook:2000} considered the $\,\star_{w}$--operation
  which
       can be defined as follows: 
       
       \centerline{$E^{\star_{w}} := \cup \{(E:H)\; |
  \;\,H \in
       \boldsymbol{f}(D) \mbox{ and } H^\star = D \}\,,
       \;\; \mbox{ for each } E \in \boldsymbol{\overline{F}}(D)\,.$}
       
       \noindent From their theory it follows that $\,\star_{w} = 
  \tilde{\star}$\,
       \cite[Corollary 2.10]{Anderson/Cook:2000}.  The relation between
       $\,\tilde{\star}\,$ and the localizing systems of ideals was
  established
       in \cite{FH}.

       \rm {\bf (b)} If $\,\Delta\,$ is a nonempty quasi-compact subset
       of $\,\mbox{Spec}(D)\,,$ \ then $\,\star_{\Delta}=
       (\star_{\Delta})_{f}\,$ and $\,{\calM}(\star_{\Delta}) =
        \Delta_{\mbox{\rm\tiny max}},$ \cite[Proposition 4.3 (B)]{FH}.

      \end{reem}

  The collection of all quasi--$\star$--ideals of a domain $\,D\,,$\
  associated to a given semistar operation $\,\star\,,$
  can be an unwieldy object.  We
  now turn to the use of ultrafilters to gain some control over this
  collection.
  A similar course was followed in \cite{CLT} for the special case
  of the $\,t$--operation.  We generalize the results given there.
  We begin with some notation/terminology/definitions.

  \begin{itemize}

  \rm \item \rm Let $\,D\,$ be a domain and let $\,{\calJ} = 
  {\calJ}(\Lambda) :=\{ J_\lambda \;|\;\,
  \lambda \in \Lambda \}\,$ be a collection of ideals of $D$.

  \rm  \item \rm Let $\,{\calU} = {\calU}(\Lambda)\,$ be an ultrafilter
  on the index set
  $\,\Lambda\,$ given above.

  \rm \item \rm For $\,I \subseteq D\,$ let $\,B(I)  := \{\lambda \in
  \Lambda \;|\;\, I \subseteq J_\lambda \}\,$.

  \rm \item \rm Let $\,J_{\calU} := \{ d \in D \;|\;\, B(d) \in 
  {\calU}
  \}\,$,\ more explicitly
  $
  \,J_{\calU} = \cup\mbox{\large \{}\ \cap\{ J_{\lambda}\;|\;\,\lambda
  \in
  B\} \;\; |  \;\,\, B \in {\calU} \ \mbox{\large \}}\,.
  $
  \ We call $\,J_{\calU}\,$ \it the $\,{\calU}$--ultrafilter
  limit of
  the collection $\,\calJ\,$.

  \end{itemize}

  %PROPOSITION 2.10
  \begin{prro}\label{pr:10} \sl Assume the
  notation/terminology/definitions given a\-bo\-ve.  Assume also that
  $\,\star\,$ is a star operation on $\,D\,$ and that each $\,J_\lambda
  \in
  {\calJ}\,$ is a $\,\star_f$--ideal [respectively, $\,\star_f$--prime]
  of
  $\,D\,$.  \ If $\,J_{\calU}\,$ is nonzero, then it is  also a 
  $\,\star_f$--ideal
  [respectively, $\,\star_f$--prime] of
  $\,D\,$.
  \end{prro} \rm

  \vskip -3pt \noindent {\bf Proof}.  The proof is the same as that 
  given in
  \cite[Proposition 2.5]{CLT} with
  the $t$--operation replaced by $\,\star\,$ an arbitrary star operation
  of finite type.
  The ``prime ideal part'' of the statement follows from \cite[Lemma
  2.4]{CLT}.
  \hfill $\Box$

  \vspace{2pt}

  %COROLLARY 2.11
  \begin{coor}\label{cor:11} \sl Generalize the setting of
  Proposition~\ref{pr:10} to the case where $\,\star\,$ is a semistar
  operation and each $\,J_{{\lambda}} \in {\calJ}\,$ is a proper
  quasi--$\star_f$--ideal [respectively, \ quasi--$\star_f$--prime] of
  $\,D\,$.
  \ If $\,J_{\calU}\,$ is nonzero, then it is also a proper 
  quasi--$\star_f$--ideal
  [respectively, \ quasi--$\star_f$--prime] of
  $\,D\,$.  \end{coor} \rm

  \vskip -3pt \noindent {\bf Proof}.  For ease of notation we set 
  $\,\boldsymbol{\ast}
  :=\star_f$ . \ As noted in Lemma~\ref{lm:1}, the
  quasi--$\boldsymbol{\ast}$--ideals of $\,D\,$ are precisely the
  contractions to $\,D\,$ of the $\,\dot{\boldsymbol{\ast}}$--ideals of
  $\,D^{\boldsymbol{\ast}}\,$ (where $ \, \dot{\boldsymbol{\ast}} \, $ 
  is
  a (semi)star operation on the domain $D^{\star_{\mbox{\tiny \it f}}}$
  defined in Remark~\ref{rk:2.1} (c)).  The result follows easily by 
  using
  Proposition~\ref{pr:10} and Lemma~\ref{lm:2}.  \hfill $\Box$

  \vspace{4pt}

  \section{Semistar Nagata Rings}

  If $\,R\,$ is a ring and $\,X\,$ an indeterminate over $\,R\,$,\ then
  the ring:

  \centerline{$
   R(X) := \{f/g \;|\;\, f,g \in R[X] \; \textrm { and } \;
  \boldsymbol{c}(g) = R\} 
  $}

  \noindent is called \it the Nagata ring of $\,R\,$ \ \rm
  \cite[Proposition 33.1]{G}.
  Some results proved in \cite[Proposition 2.1]{Kang:1989} are
  generalized in the following:

  %PROPOSITION 3.1
  \begin{prro}\label{pr:3.1}
  \sl Let $\,\star\,$ be a nontrivial semistar
  operation of an integral domain $\,D\,$. \ Set $
   N(\star) := N_D(\star) := \{h \in D[X] \; | \;\, h \not = 0  \;
  \textrm { and} \; \boldsymbol{c}(h)^\star = D^\star \} \,.
   $

  %\begin{enumerate}
  \bara

  \rm \bf \item \sl $N(\star)= N(\star_f)\,$ is a saturated 
  multiplicatively closed
  subset of
  $\,D[X]\,$.

  \rm \bf \item \sl $N(\star) = D[X] \setminus \cup \{Q[X] \; | \;\, Q 
  \in
  {\calM}(\star_f) \}\,. $

  \rm \bf \item \rm Max$(D[X]_{N(\star)}) = \{Q[X]_{N(\star)} \; |\;\, Q
  \in
  {\calM}(\star_f)\}\,.$

   \rm \bf \item \sl $D[X]_{N(\star)} = \cap\{ D_Q(X) \; |\;\, Q \in
   {\calM}(\star_f) \}\,.$

  \rm \bf \item \sl ${\calM}(\star_f)$ coincides with the canonical 
  image
  in
  \rm \, Spec$(D)$ \, \sl of the maximal spectrum of
  $\,D[X]_{N(\star)}\,;$ \ \sl i.e.
  ${\calM}(\star_{f}) = \{ M \cap D \; | \; \; M \in $ \rm
  Max$({D[X]}_{N(\star)}) \}\,.$
  %\hfill $\Box$

  %\end{enumerate}
  \eara

  \end{prro}\rm

  \vskip -3pt \noindent {\bf Proof}. (1) It is obvious that $\,N(\star) 
  =
  N(\star_f)\,$;\ the remaining part
  is a standard consequence of (2)
  \cite[Theorem 2]{Ka}.

  (2) We start by proving the following:
  \smallskip

  \noindent {\bf Claim}: \sl Let $\,h \in D[X], h \not = 0\,$.\ Then:

  \centerline{$ \boldsymbol{c}(h)^\star = D^\star \;\, \Leftrightarrow 
  \;\,
  \boldsymbol{c}(h) \not \subseteq Q\,, \; \textrm { for each } \; Q \in
  {\calM}(\star_f)\,.  $} \rm

  $(\Leftarrow)\,$ If $\,\boldsymbol{c}(h)^\star \not = D^\star ,$\ then
  $\,0 \not = \boldsymbol{c}(h) \subseteq \boldsymbol{c}(h)^\star \cap D
  \subsetneq D\,$.\  Since $\,\boldsymbol{c}(h)^{\star_f} \cap D =
  \boldsymbol{c}(h)^\star \cap D\,$ is a proper quasi--$\star$--ideal of
  $\,D\,,$\ we can find $\,Q \in {\calM}(\star_f)\,$ such that
  $\,\boldsymbol{c}(h) \subseteq \boldsymbol{c}(h)^\star \cap D 
  \subseteq
  Q\,$ (Lemma~\ref{lm:2} (1)).  This fact contradicts the assumption.

  $(\Rightarrow)\,$ is trivial.
  \smallskip

  \noindent Using the claim, we have:
  $$
  \begin {array} {ll}
   h \in N(\star) \; \Leftrightarrow \; \boldsymbol{c}(h)^\star = 
  D^\star
  &\Leftrightarrow \; \boldsymbol{c}(h) \not \subseteq Q\,, \, \textrm {
  for each} \; Q \in
  {\calM}(\star_f) \;\Leftrightarrow \\
  &\Leftrightarrow\; h \not \subseteq Q[X]\,, \, \textrm{  for each} \; 
  Q
  \in
  {\calM}(\star_f)\,.
  \end {array}
  $$

  (3) By using \cite[(4.7) and Proposition 4.8]{G}, it is sufficient to
  show
  that each prime ideal $\,H\,$ of $\,D[X]\,$ contained inside $\,\cup
  \{Q[X] \; |\;\, Q \in {\calM}(\star_f)\}\,$ is contained in
  $\,Q[X]\,$,\
  for some $\,Q \in {\calM}(\star_f)\,$. \ Let
  $\boldsymbol{c}(H)$ be the ideal generated by $ \{\boldsymbol{c}(h) 
  \;|\;\, h \in H\}\,.$
  \noindent It is easy to see that $\,\boldsymbol{c}(H)\,$ is an ideal 
  of
  $\,D\,$
  and that:
  $$
  H \subseteq \cup \{Q[X] \; | \; \, Q \in {\calM}(\star_f) \} \;
  \Rightarrow
  \; \boldsymbol{c}(H)^{\star_f} \not = D^{\star_{\mbox \tiny \it f}} =
  D^\star\,.
  $$
  \noindent As a matter of fact, if $\,\boldsymbol{c}(H)^{\star_f} =
  D^\star\,$,\ then we can find a
  polynomial $\,\ell \in \boldsymbol{c}(H)[X]\,$ such that
  $\,\boldsymbol{c}(\ell)^\star = \boldsymbol{c}(\ell)^{\star_f}=
  D^\star\,$. \ Now,
  $$
  \ell \in \boldsymbol{c}(h_1)[X] +\boldsymbol{c}(h_2)[X]+ \ldots +
  \boldsymbol{c}(h_r)[X] = (\boldsymbol{c}(h_1)+\boldsymbol{c}(h_2)+
  \ldots +\boldsymbol{c}(h_r))[X]
  $$
  \noindent with $\,(h_1, h_2, \ldots , h_r) \subseteq H\,.$ \ Since
  $\,\boldsymbol{c}(h_1)+\boldsymbol{c}(h_2)+ \ldots 
  +\boldsymbol{c}(h_r)
  \subseteq \boldsymbol{c}(H)\,$
  and $\,\boldsymbol{c}(H)\,$ is an ideal of $\,D\,,$ \ then
  $\,\boldsymbol{c}(h_1)+\boldsymbol{c}(h_2)+ \ldots 
  +\boldsymbol{c}(h_r)
  = \boldsymbol{c}(h)\,$, \ for some $\,h \in H\,.$ \ Therefore
  $\,\boldsymbol{c}(\ell) \subseteq \boldsymbol{c}(h)\,$ and thus
  $\,\boldsymbol{c}(\ell)^\star = \boldsymbol{c}(h)^\star = D^\star .$ \
  This is a contradiction, since $\,h \in H\,$ and thus
  $\,\boldsymbol{c}(h)^\star = \boldsymbol{c}(h)^{\star_f} \subseteq 
  Q\,$,
  \ for some $\,Q \in {\calM}(\star_f)\,$.
   By the fact that $\,\boldsymbol{c}(H)^{\star_f} \not =
  D^{\star_f}\,$
  we deduce that $\,\boldsymbol{c}(H) \subseteq Q\,$, for some $\,Q \in
  {\calM}(\star_f)\,$.  \ This implies that $\,H \subseteq Q[X]\,,$ for
  some $\,Q \in {\calM}(\star_f)\,$.

  (4) and (5) are easy consequences of (3), since:
  $$
  {(D[X]_{N(\star)})}_{Q[X]_{N(\star)}} = D[X]_{Q[X]} = D_{Q}(X)\,,
  $$

  \noindent (cf. also \cite[Corollary 5.3 and Proposition 33.1]{G}).
  \hfill $\Box$

  \vspace{4pt}

  We set:

  \vskip -4pt\centerline{$
  {\textstyle\rm Na}(D,\star) := D[X]_{N_{D}(\star)}
  $}

  \noindent and we call it \it the Nagata ring of $\,D\,$ with respect 
  to
  the
  semistar operation $\,\star$\,\rm.\ Obviously, $\,{\textstyle\rm
  Na}(D,\star) =
  {\textstyle\rm Na}(D,\star_f)$\,.\ If $\,\star = d\,$ is the identity
  (semi)star operation of $D$, then $
  {\textstyle\rm Na}(D,d) = D(X)\,.$

  \begin{coor} \sl Let $\,D\,$ be an integral domain, then:
   
   \centerline{$
   Q \; \mbox{is a maximal $t$--ideal of }  D \, \Leftrightarrow \,
   Q = M \cap D\,, \; \mbox{for some} \, M \in \mbox{\rm Max(Na}(D, 
  v))\,.
  $}
  \end{coor} \rm

  \vskip -3pt \noindent {\bf Proof}.  It is a straightforward 
  consequence of
  Proposition~\ref{pr:3.1} (5).  \hfill $\Box$

  %EXAMPLE3.3
  \begin{exxe}\label{ex:3.2}

  %\begin{enumerate}
  %\bara

  \rm \bf (1) \rm Let $P$ be a nonzero prime ideal of an integral domain
  $D$ and let $\,\star := \star_{\{D_P\}}\,$ be the semistar operation 
  of
  $D$
  defined as follows:

  \centerline{$ E^{\star_{\{D_P\}}} := ED_P\,, \; \,\textrm { for each} 
  \; E \in
  \boldsymbol{\overline{F}}(D)\,.
  $}

  \noindent Then, \sl it is easy to verify that:

  %\begin{enumerate}
  \balf

  \rm \bf  \item \sl${\calM}(\star_f) = \{P\}\,;$

  \rm \bf  \item \sl ${\textstyle\rm Na}(D,\star) = D_P(X)\,;$

  \rm \bf  \item \sl $\star = \star_f = \star_{sp} = \tilde{\star}\,. $

  %\end{enumerate}
  \ealf

  \rm \bf (2) \rm The previous example can be generalized as follows.
  \sl Let $\,D\,$ be an integral domain, let $\,\Delta\,$ be a nonempty
  subset
  of \rm\
  Spec$(D)$\ \sl and set $\,\star := \star_\Delta\,$. \ Let
  $\,\Delta_{\mbox{\rm \tiny max}}\,$ be the set of all the maximal
  elements of $\,\Delta\,$
  and let

  \centerline{$
  \Delta^{\downarrow} := \{ H \in \mbox{\rm Spec}(D) \; | \;\, H 
  \subseteq
  P\,, \; \textrm { for
  some } \; P \in \Delta \} \,.
  $}

  \noindent Assume that each $\,P \in \Delta\,$ is contained in some 
  $\,Q
  \in
  \Delta_{\mbox{\rm \tiny max}}\,$. \ Then, under the previous
  assumptions, $\,\star =
  \star_{\Delta_{\mbox{\rm \tiny max}}}\,$ (Lemma~\ref{lm:3} (5)) and
  moreover:

  %\begin{enumerate}
  \balf

  \rm \bf  \item \sl $\Delta \subseteq \mbox{\rm QSpec}^{\star}(D)
  \subseteq \Delta^{\downarrow}\,$, \
  thus $\,\mbox{\rm QMax}^{\star}(D) =\Delta_{\mbox{\rm \tiny max}} \,$.
  %\vspace{.13in}

  \noindent \hskip -29pt \sl Assume also that $\,\Delta_{\mbox{\rm \tiny
  max}}\,$ \ is a quasi-compact subspace
  of \ \rm Spec$(D)$.   \sl Then:

  \rm \bf \item \sl ${\textstyle\rm Na}(D,\star_\Delta) = \cap\{D_Q(X) 
  \;
  |\;\, Q
  \in \Delta_{\mbox{\rm \tiny max}} \} = \cap\{D_P(X) \; |\;\, P \in
  \Delta\}\,$.

  \rm \bf \item \sl $\widetilde{(\star_\Delta)} = \star_\Delta\,.$

  %\end{enumerate}
  \ealf
  %\end{enumerate}
  %\eara
  \end{exxe}

  \vskip -3pt \noindent {\bf Proof}. (a) If $\,P \in \Delta\,$, \ then
  $
  P^\star = P^{\star_{\Delta}} = PD_P \cap (\ \cap \{D_{P'} \;|\;\, P' 
  \in
  \Delta\,,\; P' \not \subseteq P\}\,) =
   PD_P \cap D^{\star_\Delta}\,,$
  \ and so $\, P \subseteq P^\star \cap D \subseteq PD_P \cap
  D^{\star_\Delta} \cap D = PD_P \cap D = P\,$.\
    This shows that
  $\, \Delta \subseteq \mbox{\rm QSpec}^\star(D)\,$.\ Let $\,H \in
  \mbox{\rm QSpec}^{\star}(D)\,$.\ Then $\,H^\star \not = D^\star\,$ \ 
  and
  so, for some $\,P \in
  \Delta\,$, \ $ HD_P = H^\star D_P \subseteq PD_P \not = D_P\,,
  $
  \ (cf.  also Lemma \ref{lm:2} (1)).  Henceforth, $\, H \subseteq
  HD_P
  \cap D \subseteq PD_P \cap D = P\,$.

  (b) If $\,\Delta_{\mbox{\rm \tiny max}}\,$ is quasi-compact, then
  $\, \star = \star_\Delta = \star_{\Delta_{\mbox{\rm \tiny max}}} =
  (\star_{\Delta_{\mbox{\rm \tiny max}}})_f = \star_{f}\,,$ \
  \cite[Corollary 4.6
  (2)]{FH} and so, by (a) and Lemma~\ref{lm:2} (3), $\,\Delta_{\mbox{\rm
  \tiny max}} ={\calM}((\star_\Delta)_f) = \mbox{\rm
  QMax}^{\star_{f}}(D)= \mbox{\rm QMax}^{\star}(D)\,.$\ The conclusion
  follows from Proposition~\ref{pr:3.1} (4).

  (c)  Since $\,\star_\Delta = (\star_\Delta)_f\,$ and
  $\,\Delta_{\mbox{\rm \tiny max}} = {\calM}((\star_\Delta)_f)\,$
  (Remark~\ref{rem:2.12} (b)), then we
  conclude by Corollary~\ref{cor:9} and Lemma~\ref{lm:7} (cf.  also
  Lemma~\ref{lm:3} (1)).  \hfill $\Box$

  \vspace{2pt}

  %PROPOSITION 3.4

  \begin{prro}\label{pr:3.3} \sl Let $\,\star\,$ be a nontrivial 
  semistar
  operation of an
  integral domain $\,D\,$ with quotient field $\,K\,$. \ Let
  $\,\tilde{\star} :=
  (\star_f)_{sp}\,$ be the spectral semistar operation considered in
  Corollary~\ref{cor:9}.  For each $\, E \in
  {\overline{\boldsymbol{F}}}(D)\,$,\ we have:

  \bara
  %\begin{enumerate}

  \rm \bf \item \sl  $E{\textstyle\rm Na}(D,\star) = \cap\{ED_Q(X)\; 
  |\;\,
  Q \in
  {\calM}(\star_f)\}\,$.

  \rm \bf \item \sl $E{\textstyle\rm Na}(D,\star) \cap K = \cap
  \{ED_Q \; | \;\, Q \in {\calM}(\star_f) \}\,$.

  \rm \bf \item \sl $E^{\tilde{\star}} = E{\textstyle\rm Na}(D,\star)
  \cap K$\,,\ hence if $\,E = E^\star\,,$\ then $\,E =
  E{\textstyle\rm
  Na}(D,\star) \cap K\,$.
  
  \rm \bf \item \sl  If $\,\star = \star_f\,$, then 
  $\,\tilde{\star}\, = \star_{sp}\,,$ \ hence $D^{\star_{sp}} = \cap \{ D_Q \;
  | \;\, Q \in {\calM}(\star)\}\,, $ \ and $\star_{sp}\,$ is a semistar
  operation of finite type.

  %\end{enumerate}
  \eara

  \end{prro} \rm

  \vskip -3pt \noindent {\bf Proof}.  (1) By Proposition~\ref{pr:3.1} 
  (3) and
  \cite[Chapitre 2, Corollaire 3, p.  112]{B} (see also the proof of
  Proposition~\ref{pr:3.1} (4)), we have:
  $$
  \begin{array}{ll}
  E{\textstyle\rm Na}(D,\star) &= \; \cap \{(ED[X]_{N(\star)})_M \; 
  |\;\,
  M
  \in \mbox{\rm Max}(D[X]_{N(\star)} \} = \\
  &= \; \cap \{(ED[X]_{N(\star)})_{Q[X]_{N(\star)}} \;|\;\, Q \in
  {\calM}(\star_f) \} = \\
  & = \cap  \; \{ED[X]_{Q[X]} \,| \,\, Q \in {\calM}(\star_f)\} = \cap
   \{ ED_Q(X)
  \,| \,\, Q \in {\calM}(\star_f)\}\,.
  \end{array}
  $$

  (2) By using (1) and \cite[Proposition 33.1 (4)]{G}, we have:
  $$ \begin{array}{rl}
  ED[X]_{N(\star)} \cap K \hskip -4.pt &= \cap  \{ ED_Q(X) \; | \; \; Q
  \in
  {\calM}(\star_f) \} \ \cap \ K = \\
  \hskip -4.0pt &=  \cap  \{ED_Q(X) \cap K \; | \; \; Q \in 
  {\calM}(\star_f) \} =\\
  \hskip -4.0pt &= \cap  \{ED_Q \; | \; \; Q \in {\calM}(\star_f)\}.
  \end{array}
  $$

  (3) From Corollary~\ref{cor:9} (2) we know that $
  E^{\tilde{\star}} = \cap \{ ED_Q \; | \;\, Q \in {\calM}(\star_f) \}
  \,,
  $
  thus the first statement of (3) is a straightforward
  consequence of (2).
  Since $\,\tilde{\star} \leq \star_f \leq \star\,$
  (Corollary~\ref{cor:9} (1)), then obviously if $\,E = E^\star\,$,\ 
  then
  $\,E = E^{\tilde{\star}}\,$ and thus $\,E = E{\textstyle\rm 
  Na}(D,\star)
  \cap K\,$.

   For a direct proof of the second statement of (3), it is
  enough to
  show that, if $\,z \in E{\textstyle\rm Na}(D,\star) \cap K\,,$ \ then
  $\,z \in
  E\,$. \ Let $\, g, h \in D[X]\, $ with $\, h \not = 0\,,$ $\,
  \boldsymbol{c}(h)^\star = D^\star \,$ and
  $\,zh = g\,$.  Then $
  \boldsymbol{c}(g)^\star = \boldsymbol{c}(zh)^\star =
  z\boldsymbol{c}(h)^\star =
  zD^\star = (zD)^\star \subseteq E^\star = E \,.$
  
  (4) follows directely from the definitions, \cite[Proposition
  3.2]{FH} and from (2) and (3).  \hfill $\Box$

%   %\vspace{4pt}
% 
% 
% 
% 
% 
% 
% 
% 
% 
% 
% 
% 
% 
% 
% 
% 
%   %COROLLARY 3.5
%   \begin{coor}\label{cor:3.4} \sl Let $\,\star\,$ be a nontrivial 
%   semistar
%   operation of an integral domain $\,D\,$ with quotient field $\,K\,$.\
%   We
%   can always define the map $\,\tilde{\star}\, :
%   \boldsymbol{\overline{F}}(E) \rightarrow 
%   \boldsymbol{\overline{F}}(E)\,$
%   by:
% 
%   \centerline{$
%   E^{\,\tilde{\star}\,} := E{\textstyle\rm Na}(D,\star) \cap 
%   K,
%   \; \textrm {  for each }
%   E \; \in \boldsymbol{\overline{F}}(D)\,.
%   $}
% 
%   \noindent If $\,\star = \star_f\,$, then:
% 
%   %\begin{enumerate}
%   \bara
% 
%   \rm \bf \item \sl  $\,\tilde{\star}\,$ is a spectral
%   semistar operation of $\,D\,$ with
%   $\,\tilde{\star}\, = \star_{sp}\,.$
% 
%   \rm \bf \item \sl $D^{\star_{sp}} = \cap \{ D_Q \; | \;\, Q \in 
%   {\calM}(\star)\}\,.$
% 
%   \rm \bf \item \sl $\star_{sp}\,$ is a semistar operation of finite 
%   type.
% 
%   %\end{enumerate}
%   \eara
% 
% %   \end{coor}
%   \rm
% 
% 
% 
% 
% 
% 
% 
% 
% 
% 
% 
%   \vskip -3pt \noindent {\bf Proof}.  (1) is a restatement of 
%   Proposition~\ref{pr:3.3}
%   (3).
% 
%   (2) follows from (1), since
%   $\, D^{\star_{sp}} = D^{\tilde{\star}\,} = {\textstyle\rm
%   Na}(D,\star) \cap
%   K\,$.
% 
%   (3) Note that the canonical image of \ Max$({\textstyle\rm
%   Na}(D,\star))\,$
%   in \ Spec$(D)\,$ coincides with $\,{\calM}(\star)\,$
%   (Proposition~\ref{pr:3.1} (5)).  Since\ Max$({\textstyle\rm
%   Na}(D,\star))\,$ is compact then $\,{\calM}(\star)\,$ is a compact
%   subspace of\ Spec$(D)\,$;\ we conclude by \cite[Corollary
%   4.6 (2)]{FH}.  \hfill $\Box$
% 

  \vspace{3pt}

%   Note that the assumption $\,\star = \star_f\,$ in
%   Corollary~\ref{cor:3.4} is not
%   really restrictive, since it is obvious that $\,{\textstyle\rm
%   Na}(D,\star) =
%   {\textstyle\rm Na}(D,\star_f)\,$.

  %  NEW COROLLARY 3.5
  \begin{coor}\label{cor:3.5} \sl  Let $\,\star\,$ be a nontrivial
  semistar
  operation of an integral domain $\,D\,,$\ let $\,\tilde{\star}\,$ be 
  the
  semistar operation of $\,D\,$ considered in Corollary~\ref{cor:9},
  and let $\,\dot{\tilde{\star}}\ (:=
  \dot{\tilde{\star}}^{D^{{\tilde{\star}}}}\ )\,$ be the (semi)star
  operation of $\,D^{{\tilde{\star}}}\,$ associated to
  $\,{\tilde{\star}}\,$ (Remark~\ref{rk:2.1}). 

  %\begin{enumerate}
  \bara

  \rm \bf \item \sl $(\tilde{\star})_f = \tilde{\star} =
  (\tilde{\star})_{sp}= \tilde{\tilde{\star}}\,$.

  \rm \bf \item \sl ${\calM}(\star_f) = {\calM}(\tilde{\star})\,$.

  \rm \bf \item \sl ${\textstyle\rm Na}(D,\star) = {\textstyle\rm
  Na}(D,\tilde{\star}) = {\textstyle\rm
  Na}(D^{\tilde{\star}},\dot{{\tilde{\star}}}) \,$.

  %\end{enumerate}
  \eara

  \end{coor} \rm

  \vskip -3pt \noindent {\bf Proof}.  (1) follows easily from
  Proposition~\ref{pr:3.3} (4) and Lemma~\ref{lm:7}.

  (2) Let $\, Q \in {\calM}(\star_{f})\,$.\ Then $\, Q^{\star_f} \cap D 
  =
  Q\,$.\ Since $\,\tilde{\star} \leq \star_f\,$ (Corollary~\ref{cor:9}
  (1)), then necessarily $\,Q^{\tilde{\star}} \cap D = Q\,$.\ By (1) we
  know that $\,\tilde{\star}\,$ is a semistar operation of finite type.
  Hence we know that the quasi--$\tilde{\star}$--ideal $\,Q\,$ of 
  $\,D\,$
  is contained in some $\,H \in {\calM}(\tilde{\star})\,$
  (Lemma~\ref{lm:2} (1)).

  Conversely, let $\,H \in {\calM}(\tilde{\star})\,$. \ Then
  $\,H =
  H^{\tilde{\star}} \cap D = \cap \{ HD_Q \; | \; \, Q \in 
  {\calM}(\star_f)\} \cap D\,$ (Corollary~\ref{cor:9} (2)).  In 
  particular,
  we have $\,HD_Q \not = D_Q\,$ for some $\,Q \in {\calM}(\star_f)\,,$ \
  since otherwise $\,H^{\tilde{\star}} \cap D\,$ would be equal to
  $\,D\,$. \ Therefore, $\,H = H^{\tilde{\star}} \cap D \subseteq HD_Q
  \cap D \subseteq QD_Q \cap D = Q\,$, \
   for some $\,Q \in {\calM}(\star_f)\,$.

  \noindent By the previous properties, we deduce immediately that
  $\,{\calM}(\star_f) = {\calM}(\tilde{\star})\,$.

  (3) Since, from (1), we know that $\,\tilde{\star}\,$ is a semistar
  operation of finite type, then, by Proposition~\ref{pr:3.1} (4),
  $
  {\textstyle\rm Na}(D,\tilde{\star}) = \cap \{D_H(X) \; | \;\, H  \in
  {\calM}(\tilde{\star})\}\,.
  $
  \ Since, by (2), we know that $\,{\calM}(\tilde{\star}) =
  {\calM}(\star_f)\,$, \ then we have $\,{\textstyle\rm
  Na}(D,\tilde{\star}) = {\textstyle\rm Na}(D,\star_f) = {\textstyle\rm
  Na}(D,\star)\,$.

  \vskip 2pt

  \noindent \bf Claim. \rm $\,{\calM}(\dot{\tilde{\star}}) = \{ QD_{Q}
  \cap D^{\tilde{\star}} \; | \;\, Q \in {\calM}({\star}_{f}) \}\,.$
  \vskip 2pt

  If $\,M \in {\calM}(\dot{\tilde{\star}})\,,$ \ then $\, M =
  M^{\dot{\tilde{\star}}}= M^{{\tilde{\star}}} = \cap \{ MD_{Q} \; | 
  \;\,
  Q \in {\calM}(\tilde{{\star}}) = {\calM}({\star}_{f}) \}\,,$ \
  hence $\,MD_{Q} \neq D_{Q}\,,$ \ thus  $\,M \subseteq QD_{Q} \cap
  D^{\tilde{\star}}\,,$ \ for some $\,Q \in  {\calM}({\star}_{f})\,.$ \
  On the other hand, it is easy to verify that $\,QD_{Q} \cap
  D^{\tilde{\star}}\,$ is a $\,\dot{{\tilde{\star}}}$--ideal of
  $\,D^{\tilde{\star}}\,$ (Lemma~\ref{lm:1}),  hence the claim is 
  proved.

  The last equality in (3) is a straightforward consequence
  of Proposition~\ref{pr:3.1} (4), of  the Claim and of the fact that
  $\, D^{\tilde{\star}}_{QD_{Q} \cap
  D^{\tilde{\star}}} = D_{Q}\,,$ \ for each $\,Q \in  
  {\calM}({\star}_{f})\,.$
  \hfill $\Box$

  %REMARK 3.6
  \begin{reem} \label{new} \rm Let $\,\star\,$ be a nontrivial semistar
  operation
  of an integral domain $\,D\,,$ \ let $\,\tilde{\star}\,$ be the 
  semistar
  operation considered in Corollary~\ref{cor:9} and let $\,\dot{\star} 
  :=
  \dot{\star}^{D^{\star}}$\ [respectively, $\,{\tilde{\dot{\star}}} :=
  {\tilde{\dot{\star}}}^{D^{\star}}]\ $\ be the (semi)star operation of
  $\,D^{\star}\,$
  associated to $\,\star\,$
  [respectively, $\,\dot{\star}\ $] and defined in Remark~\ref{rk:2.1}
  (c)
  [respectively,   Corollary~\ref{cor:9}].\ Then, in general, \sl the
  semistar Nagata ring $\,{\textstyle\rm
  Na}(D,\star) = {\textstyle\rm Na}(D,\tilde{\star})\,$ is different 
  from
  the star Nagata ring $\,{\textstyle \rm Na}(D^{\star},\dot{\star}) =
  {\textstyle \rm Na}(D^{\star},\tilde{\dot{\star}})\,.$ \rm

  \smallskip

  Let $\,L\,$ be a field and $\,X,\ Y, \ U,\,$ and $\, Z\,$ 
  indeterminates
  over
  $\,L\,.$ \ Set:

  \centerline{$ V:= L(X)[Y]_{(Y)} = L(X) + M\,, \;\, M:= 
  YL(X)[Y]_{(Y)}\,,
  \;\, R:= L + M\,.
  $}

   It is well known that $\,R\,$ and $\,V\,$ are 1-dimensional
  local domains with the same field of quotients $\,F:= L(X,Y)\,$, \ 
  with
  the
  same set of prime ideals $\,\{(0), M\}\,,$ \ and moreover $\,V\,$ is a
  discrete valuation domain, $\,\mbox{dim}(V[U])=2\,$ and
  $\,\mbox{dim}(R[U])=3\,.$ \ The last property follows from the fact 
  that
  in $\,R[U]\,$ we have the following inclusions of prime ideals:
  $$
  \begin{array}{ll}
  (0) \subset Q_{1} \hskip -5pt &:= (U-X)F[U] \cap R[U] = (YU-YX)F[U] 
  \cap
  R[U] \subset\\
  \hskip -5pt & \subset Q_{2}:= M[U] \subset Q_{3}:=(M,U)R[U] \,,
  \end{array}
  $$
  \noindent but, a similar inclusion does not hold in $\,V[U]\,$:

  \centerline{$
  %\begin{array}{ll}
  (0) \subset P_{1} := (U-X)F[U] \cap V[U] = (U-X)V[U] \not\subseteq  
  M[U]
  \,;
  %\\
  %\hskip -5pt &\not\subseteq
  %\end{array}
  $}

  \noindent in fact, more generally, no height 1 prime ideal $\,P\,$ of
  $\,V[U]\,,$
  \ with $\,P \cap V = (0)\,,$ is
  contained in $\,M[U]\,$ \cite[Theorem 39, Theorem 68; Exercise 18, 
  page
  42]{Ka}.

   Set $\, D:= R(U)\,,$ \ $\, T:= V(U)\,,$ \ and let $\,\star:=
  \star_{\{T\}}\,$ be the semistar operation of $\,D\,$ considered in
  Remark ~\ref{rk:2.1} (d).  Since all the prime ideals of $\,D\,$
  [respectively, of $\,T\,$], different from $\,(0)\,$ and $\,M[U]\,,$ \
  are
  of the type $\,fF[U] \cap D\,$ [respectively, $\,fF[U] \cap T\,$], 
  where
  $\,f \in F[U]\,$ is irreducible \cite[Theorem 36]{Ka}, then it follows
  that the canonical map \ Spec$(T) \rightarrow$ \ Spec$(D)$ is a
  bijection.  From this fact we deduce immediately that 
  $\,{\calM}(\star)
  = \mbox{Max}(D)\,,$ \ and, hence, that $\,\tilde{\star} = d_{D}\,,$ \
  where $\,d_{D}\,$ is the identity (semi)star operation of $\,D\,.$

    Note that, in the present situation, $\,
   D^{\star}=T\,$
  and it is obvious that
   $\,\dot{\star}\,$ coincides
  with $\,d_{T}\,,$ \ where $\,d_{T}\,$
    is  the identity (semi)star operation of $\,T\,,$ \ and hence
    $\,\dot{\tilde{\star}}= \dot{d}_{D} = d_{T}\,.$ \
  \ We deduce that $\, {\calM}(\dot{\star}) = 
  {\calM}(\dot{\tilde{\star}})=\mbox{Max}(T) \,$ and, obviously, that
  $\,\widetilde{d_{T}}= d_{T}.$

  Furthermore, note that $\, P_{1}T\,$ is a maximal ideal of $\,T\,$
   (because
  $\,P_{1} \not\subseteq M[U] \subset V[U]\,$), \ but $\,P_{1}T \cap D =
  Q_{1}D\,$ is not maximal in $\,D\,$ (because $\,Q_{1} \subset M[U]
  \subset R[U]\,).$ \ Therefore the  statement in 
  Lemma~\ref{lm:2} (3) is not reversible.  Note also that this example 
  shows
  that $\,(Q_{1}D)^{\star} = Q_{1}T \subsetneq P_{1}T\,$ (because
  $\,Q_{1}V[U] = (YU-YX)V[U] \subsetneq P_{1}=(U-X)V[U]\,$) and so
  $\,(Q_{1}D)^{\star}\,$ is not a $\,\dot{\star}$--prime of
  $\,D^{\star}=T\,,$ \ even though $\,(Q_{1}D)^{\star} \cap D = Q_{1}D
  \,$ is a quasi--$\star$--prime of $\,D\,;$ \ so also the 
  statement in Lemma~\ref{lm:2} (4)  is not reversible.

   Finally note that:
  $$
  \begin{array}{ll}
  {\textstyle\rm Na}(D,\star) &=
   {\textstyle\rm Na}(D,\tilde{\star}) = {\textstyle\rm
  Na}(D^{\tilde{\star}},\dot{\tilde{\star}})=
    {\textstyle\rm Na}(D,d_{D})
   =D(Z) \subsetneq \\
   &\subsetneq
   {\textstyle\rm Na}(D^{\star},\dot{\star}) =
   {\textstyle\rm Na}(D^{\star},\tilde{\dot{\star}})=
   {\textstyle\rm Na}(T,d_{T}) = T(Z)\,.
   \end{array}
   $$

  \end{reem}

  %COROLLARY 3.7
  \begin{coor}\label{cor:3.6} \sl Let $\,\star\,$ be a semistar 
  operation
  of
  an integral domain $\,D\,$. \  Assume that $\,\Pi^\star \not =
  \emptyset\,$
  and that $\,\star\, $ is quasi-spectral.  Then:

  \centerline{$
  {\textstyle\rm Na}(D,\star) = {\textstyle\rm Na}(D,\star_{sp}) =
  {\textstyle\rm Na}(D,\tilde{\star})\,.
  $}
  \end{coor} \rm

  \vskip -3pt \noindent {\bf Proof}. Under the present assumptions, we 
  can define the
  nontrivial semistar operation $\,\star_{sp}\,$ and we have that
  $
  \tilde{\star} = (\star_f)_{sp} \leq \star_{sp} \leq \star \,,
  $
  (Lemma~\ref {lm:7} (2) and Corollary~\ref{cor:9}).  Since it is
  easy to
  see that $
  \star_1 \leq \star_2 \,$ implies that ${\textstyle\rm Na}(D,\star_1)
  \subseteq {\textstyle\rm Na}(D,\star_2)\,,
  $
  \ then the conclusion follows immediately from
  Corollary~\ref{cor:3.5} (3).
  \hfill $\Box$

  \vspace{4pt}

  The content of Proposition~\ref{pr:3.1} (5) is that, when the maximal
  ideals of
  \,Na$(D,\star)\,$ are contracted to $\,D\,$, \ the result is exactly 
  the
  prime
  ideals of $\,D\,$ in $\,{\calM}(\star_f)\,$. \ We now prove that this
  result can be
  reversed: the maximal ideals of \,Na$(D,\star)\,$ can be obtained by
  extending to \,Na$(D,\star)\,$ the prime ideals of $\,D\,$ in 
  $\,{\calM}(\star_f)\,$.  In particular:

  %THEOREM 3.8
  \begin{thee}\label{tm:3.7}  \sl Let $\,\star\,$ be a semistar 
  operation
  of
  an integral domain $\,D\,$. \ Then $\,\mbox{\rm Max(Na}(D,\star)) =
  \{QD_Q(X)
  \cap \mbox{\rm Na}(D,\star) \;| \;\, Q \in {\calM}(\star_f)\}\,.$
  \end{thee} \rm

  \vskip -3pt \noindent {\bf Proof}.  Proposition~\ref{pr:3.1} (3) 
  indicates that the
  maximal
  ideals of \,Na$(D,\star)\,$ are exactly the ideals of
  the set $\,\{Q[X]_{N(\star)}\; | \;\, Q \in {\calM}(\star_f) \}\,$.\
  The
  result follows easily since these ideals are maximal in\ Na$(D, 
  \star)\
  $and
  are each contained in an ideal of the form $\,QD_Q(X)\,$, where $\,Q 
  \in
  {\calM}(\star_f)\,$.  \hfill $\Box$

  \vspace{4pt}

  Note that the previous result indicates a strong similarity between 
  the
  Nagata
  rings and the Kronecker function rings associated to a given 
  semi\-star
  operation.  In particular, the maximal spectrum of each ring consists 
  of
  restrictions of the maximal ideals of local overrings of the form
  $\,R(X)\,$
  where $\,R\,$ is a local overring of $\,D\,$ (cf.  
  Theorem~\ref{tm:3.7}
  and
  \cite[Theorem 3.5]{FL1}).  The difference is that, in the Kronecker
  case, the overrings $\,R\,$ are valuation overrings of $\,D\,$ and, in
  the Nagata case, they are localizations of $\,D\,$ at certain prime
  ideals.  This is perhaps an indication that the Nagata and Kronecker
  constructions are actually each special cases of a more general
  construction involving more general classes of overrings.

  \vspace{4pt}

  We now turn our attention to the question of valuation overrings.  The
  notion that we recall next is due to P. Jaffard
  \cite{J} (cf.  also \cite{HK1}, \cite{HK3}, \cite{FL2}).
  For a domain $\,D\,$ and a semistar operation $\,\star\,$ on $\,D\,$, 
  \
  we say that a valuation overring $\,V\,$ of $\,D\,$ is \it a
  $\,\star$--valuation overring of $\,D\,$ \rm provided $\,F^\star
  \subseteq FV\,,$ \ for each $\,F \in \boldsymbol{f}(D)\,$.\ Note that,
  by definition the $\,\star$--valuation overrings coincide with the
  $\,\star_{f}$--valuation overrings; by \cite[Proposition 3.3]{FL2} the
  $\,\star$--valuation overrings also coincide with the
  $\,\star_{a}$--valuation overrings.

  %THEOREM 3.9
  \begin{thee}\label{tm:3.8} \sl Let $\,D\,$ be a domain and let
  $\,\star\,$ be a
  semistar operation on $\,D\,$.  \ A  valuation overring $\,V\,$ of
  $\,D\,$ is a $\,\tilde{\star}$--valuation overring of $\,D\,$ if and
  only if $\,V\,$ is an overring of $\,D_P\,,$ for some $\,P\in
  {\calM}(\star_f)\,$.  \end{thee} \rm

  \vskip -3pt \noindent {\bf Proof}.  To avoid the trivial case, we 
  assume that $\,V
  \neq K\,.$ \ First suppose that $\,V\,$ is a valuation overring of
  $\,D_P\,,$ for
  some $\,P \in {\calM}(\star_f)\,.$ \ It is clear from the definition 
  of
  $\,\tilde{\star}\,$ that $\,V\,$ is a $\,\tilde{\star}$--valuation
  overring of $\,D\,$.

  Now assume that $\,V\,$ is a  $\,\tilde{\star}$--valuation overring of
  $\,D\,$.  \ Let $\,M\,$ be the maximal ideal of $\,V\,$ and let $\,P 
  :=
  M \cap D\,$.  \ We need to show that $\,P\,$ is contained in a prime
  $\,Q \in {\calM}(\star_f)\,$.\ We consider two cases.

  \noindent \bf Case 1. \rm \, Suppose that there is a finitely 
  generated
  ideal
  $\,J\,$ of $\,D\,$ contained in $\,P\,$ such that $\,J \not \subseteq
  Q\,,$ for each $\,Q \in {\calM}(\star_f)\,$.\ Then $
  \,J^{\tilde{\star}} = \cap \{ JD_Q \;| \;\, Q \in {\calM}(\star_f) \}
  = \cap \{ D_Q  \;| \; \, Q \in {\calM}(\star_f)\} =
  D^{\tilde{\star}}\,.$
  \  However, $\,JV \subseteq PV\,$ is a proper ideal of $\,V\,$
  and so
  cannot contain $\,J^{\tilde{\star}}\,$.  \ This contradicts our
  assumption that $\,V\,$ was a $\,\tilde{\star}$--valuation overring of
  $\,D\,$.  \ We conclude that no such ideal $\,J\,$ can exist.

  \noindent \bf Case 2. \rm \, Suppose that every finitely generated 
  ideal
  $\,J\,$
  which is contained in $\,P\,$ is also contained in some ideal $\,Q \in
  {\calM}(\star_f)\,$.\ Let $\, {\calM}(\star_f) = \{Q_\lambda\; | \; \,
  \lambda
  \in \Lambda \}\,$.\ Note that, by assumption, for any finitely 
  generated
  ideal $\,J \subseteq P\,$, the set $\,B(J) := \{\lambda \in
  \Lambda \;|\;\, J\subseteq Q_{\lambda}\}\,$ is not empty.  Let 
  ${\calU}\,$ be an ultrafilter on $\,\Lambda\,$ which contains each set
  $\,B(J)\,,$\ where $\,J\,$ runs through the finitely generated ideals
  contained in $\,P\,.$ \ Such an ultrafilter exists because the
  intersection of any finite collection of sets $\,\{ B(J_1), B(J_2),
  \ldots , B(J_n) \}\,$ is simply $\,B(J_1 + J_2 + \ldots + J_n)\,$ and 
  so
  must be nonempty.  Then the ultrafilter limit ideal $\,Q_{{\calU}}\,$
  of the collection of prime ideals $\,{\calM}(\star_f)\,$ must be a
  $\,\star_f\,$--prime (Proposition~\ref{pr:10}) and it also clearly
  contains $\,P\,$.\ This completes the proof.  \hfill $\Box$

\newpage
\section{Semistar Kronecker function rings}

  Let $\,\star\,$ be a semistar operation on an integral domain 
  $\,D\,$.\
  We say that $\,\star\,$ is \it an e.a.b. (endlich arithmetisch
  brauchbar)
  semistar operation of $\,D\,$ \rm if, for all $\, E,F,G \in
  \boldsymbol{f}(D)\,,$ \ 
  $(EF)^\star \subseteq (EG)^\star$\ implies that\ $F^\star \subseteq
  G^\star\,,$ \ \cite[Definition 2.3 and Lemma 2.7]{FL1}.

  It is possible to associate to any semistar operation $\,\star\,$ of
  $\,D\,$ an e.a.b.
  semistar operation of finite type
  $\, \star_a\, $ of $\,D\,$,\ defined as
  follows:
  $$
  \begin {array} {ll}
  F^{\star_a} &:= \cup\{((FH)^\star:H^\star) \; \ | \; \, \; H \in
  \boldsymbol{f}(D)\}\,, \; \textrm { for each} \; F \in
  \boldsymbol{f}(D)\,;\\
  E^{\star_a} &:= \cup\{F^{\star_a} \; | \; \, F \subseteq E\,,\; F \in
  \boldsymbol{f}(D)\}\,, \; \textrm { for each} \; E \in
  {\overline{\boldsymbol{F}}}(D)\,.
  \end{array}
  $$
  \noindent The semistar operation $\,\star_a\,$ is called \it the 
  e.a.b.
  semistar
  operation associated to $\,\star\,$ \rm \cite[Definition 4.4]{FH}.  
  Note
   the previous construction is essentially due to P. Jaffard \cite{J}
  (cf.
   also F. Halter-Koch \cite{HK2}).  Note that $\,D^{\star_a}\,$ is
   integrally closed and contains the integral closure of $\,D\,$ in
   $\,K\,$ \cite[Proposition 4.5]{FL1} (cf. also \cite{HK1} \cite{HK3},
   \cite{OM1} and \cite{FL2}). When $\,\star = v\,$, then 
  $\,D^{v_{a}}\,$
  coincides
   with the pseudo-integral closure of $\,D\,$ introduced by D.F.
   Anderson, Houston and Zafrullah \cite{AHZ}.

  \vspace{4pt}

  If $\,\star\,$ is a semistar operation of an integral domain $\,D\,$,\
  then we
  call \it the Kronecker function ring of $\,D\,$ with respect to
  $\,\star\,$ \rm the
  following domain:
  $$
  \begin {array} {rl}
  \mbox{Kr}(D,\star) := \{ f/g \;  \,|\, & f,g \in D[X] \setminus \{0\}
  \;\;
  \mbox{ \rm and there exists } \; h \in D[X] \setminus \{0\} \;   \\ &
    \mbox{ \rm such that } \; (\boldsymbol{c}(f)\boldsymbol{c}(h))^\star
    \subseteq (\boldsymbol{c}(g)\boldsymbol{c}(h))^\star \,\} \, \cup\,
    \{0\}\} \,,
  \end{array}
  $$
  \noindent \cite[Theorem 5.1]{FL1} (cf. also \cite{HK3}, \cite{M} and
  \cite{OM3}).

  \vspace{4pt}

  In the following statement we collect some of the main properties
  related to the Kronecker
  function ring of an integral domain with respect
  to a semistar operation (cf.  \cite[Proposition 3.3, Theorem 3.11, 
  Proposition 4.5, Theorem 5.1 and the proof of Corollary 5.2]{FL1}).

  %PROPOSITION 4.1

  \begin{prro}\label{pr:4.1} \sl Let $\,\star\,$ be a semistar operation
  of an
  integral domain $\,D\,$ with quotient field $\,K\,,$ \ let 
  $\,\star_a\,$
  be the
  e.a.b. semistar operation of $\,D\,$ associated to $\,\star\,$ and let
  $\,{\dot{{\star}}_{a}} \ (= \ {\dot{{\star}}_{a}}^{D^{\star_{a}}} \ 
  )\,$
  be the
  (semi)star operation of $\,D^{\star_{a}}\,$ associated
  to $\,{\star_{a}}\,$ and defined in Remark~\ref{rk:2.1} (c).  Then:

  %\begin{enumerate}
  \bara
  \rm \bf \item \sl $\star_f \leq \star_a\,.$

  \rm \bf \item \sl ${\textstyle\rm Kr}(D,\star) = {\textstyle\rm
  Kr}(D,\star_f) = {\textstyle\rm Kr}(D,\star_a) = {\textstyle\rm
  Kr}(D^{\star_{a}},{\dot{{\star}}_{a}} )\,.$

  \rm \bf \item \sl ${\textstyle\rm Kr}(D,\star)$ is a B\'ezout domain
  with
  quotient field $\,K(X)\,$.

  \rm \bf \item \sl ${\textstyle\rm Na}(D,\star) \subseteq 
  {\textstyle\rm
  Kr}(D,\star)\,.$

  \rm \bf \item \sl $E^{\star_a} = E{\textstyle\rm Kr}(D,\star) \cap
  K\,,$ \, for each $\,E \in \boldsymbol{\overline{F}}(D)\,.$ \hfill 
  $\Box$
  %\end{enumerate}
  \eara

  \end{prro}\rm

  %REMARK 4.2
  \begin{reem} \label{rk:4.2} \rm Note that if $\,\star\,$ is not the
  trivial semistar
  operation, then $\,\star_a\,$ is also different from the trivial
  semistar
  operation.  As a matter of fact, if $\,D^\star \not = K\,$,\, then
  $\,D^{\star_a} = \cup \{ (H^\star: H^\star)\; |\;\, H \in
  \boldsymbol{f}(D) \}
  \not = K\,$.\,  Otherwise $\,(H^\star: H^\star)\,$ would be equal to
  $\,K\,$,\ for
  some $\,H \in \boldsymbol{f}(D)\,$ and $\,H \subseteq D\,$. \ This
  implies easily
  that $\,H^\star = K\,$ and this contradicts the assumption that
  $\,D^\star
  \not = K\,$.
  \end{reem}

  %THEOREM 4.3
  \begin{thee} \label{thm:4.3} \sl Let $\,\star\,$ be a nontrivial
  semistar operation
  of an integral domain $\,D\,$.\
    Assume that $\,\star = \star_f\,$.\ Let $\,\star_a\,$ be the
  e.a.b. semistar operation of finite type canonically associated to
  $\,\star\,$.

  %\begin{enumerate}
  {\bara

  \rm \bf \item \sl Let $\,(W,N)\,$ be a nontrivial valuation overring 
  of
  $\,{\textstyle\rm Kr}(D,\star)\,$.\ Set $\,N_{0} := N \cap D\,$ and 
  let
  $\,N_{1}:= N \cap D[X]\,$.\ Then:

  %begin{enumerate}
  {\balf
  \rm \bf \item \sl $N_{1} = N_{0}[X]\,,$ \ $\,N \cap {\textstyle\rm
  Na}(D,\star) = N_{0}{\textstyle\rm Na}(D,\star) = N_{1}{\textstyle\rm
  Na}(D,\star)\,$ and $\,N \cap {\textstyle\rm
  Na}(D,\star_{a}) = N_{0}{\textstyle\rm Na}(D,\star_{a}) =
  N_{1}{\textstyle\rm
  Na}(D,\star_{a})\,$

  \rm \bf \item \sl $N_{0}\,$ is a quasi--$\star_a$--prime ideal
  (in particular, a quasi--$\star$--prime ideal) of $\,D\,.$

  %\end{enumerate}
  \ealf}

  \rm \bf \item \sl If $\,P\,$ is a quasi--$\star_a$--prime ideal of
  $\,D\,,$\
    then there exists a quasi--$\star_a$--maximal ideal $\,Q\,$ of 
  $\,D\,$
  and a valuation overring
  $\,(W,N)\,$ of $\,{\textstyle\rm Kr}(D,\star)\,$ such that $ \,P
  \subseteq Q = N \cap D\,.$

  \rm \bf \item \sl ${\calM}(\star_a)\,$ is contained in the canonical
  image
  in $\,D\,$ of\, \rm Max$({\textstyle\rm Kr}(D,\star))\,$.

  \rm \bf \item \sl For each $\,Q \in {\calM}(\star_a)\,,$\ there exists
  a $\,\star$--valuation overring $\,(V,M)\,$ of $\,D\,$ dominating
  $\,D_Q\,.$
  %and $\, Q{\textstyle\rm Kr}(D,\star) = MV(X) \cap {\textstyle\rm
  % Kr}(D,\star)\,.$

  %\end{enumerate}
  \eara}

  \end{thee} \rm

  \vskip -3pt \noindent {\bf Proof}.  As usual, we denote by $\,K\,$ 
  the field of
  fractions of
  $\,D\,$. \ It is obvious that:
  (1, a) $\, N_{0}[X] \subseteq N_{1} = N \cap D[X]\,$, \ and if
  $\,f := f_0 + f_1X + \ldots + f_rX^r \in N \cap D[X]\,$,\ then $
  {\boldsymbol c}(f){\textstyle\rm Kr}(D,\star) = f{\textstyle\rm
  Kr}(D,\star)
  \subseteq N\,
  $  \cite[Theorem 3.11 (2) and Theorem 5.1 (2)]{FL1}.
  Therefore,
  $\,f_i \in N \cap D = N_{0}\,,$\ for each $\,i\,$ with $\, 0 \leq i 
  \leq
  r\,$. \ This fact implies that $\,f \in N_{0}[X]\,$.

  \noindent Since $\,{\textstyle\rm Na}(D,\star)\,$ (and 
  $\,{\textstyle\rm
  Na}(D,\star_a)$\,) is a ring of fractions of $\,D[X]\,$
  and $\,N_{0}[X] = N_{1} = N \cap D[X]\,,$\, we have immediately that \
  %$$\begin{array}{ll}
  %&
  $N_{0}{\textstyle\rm Na}(D,\star)=N_{1}{\textstyle\rm Na}(D,\star) = N
  \cap {\textstyle\rm
  Na}(D,\star) \;\; \textrm{ and }$ 
  %&
  $ N_{0}{\textstyle\rm Na}(D,\star_a)=N_{1}{\textstyle\rm 
  Na}(D,\star_a) =
  N \cap {\textstyle\rm Na}(D,\star_a)\,.$
  %\end{array} $$

  (1, b)  Recall that $\,{\textstyle\rm Na}(D,\star_a) \subseteq
  {\textstyle\rm Kr}(D,\star_a) =$
  $ {\textstyle\rm Kr}(D,\star)\,$ and  $\,
   N_{0}^{\star_a} = N_{0}{\textstyle\rm Kr}(D,\star) \cap K\,,
   $
  \ (Proposition~\ref{pr:4.1} (2), (4) and (5)).  Since
  $\,N_{0}{\textstyle\rm Kr}(D,\star) \subseteq N \cap {\textstyle\rm
  Kr}(D,\star)\,$, \ then
  $\,
  N_{0}^{\star_a} \subseteq  N \cap {\textstyle\rm Kr}(D,\star) \cap K =
  N \cap D^{\star_a}\,,
  $
  \  (Proposition~\ref{pr:4.1} (5)).  Thus
  $\,N_{0}^{\star_a} \cap D \subseteq N \cap D^{\star_a} \cap D = N \cap
  D = N_{0}\,$, \, i.e. $\,N_{0}\,$ is a quasi--$\star_a$--prime ideal 
  of
  $\,D\,.$ \ Since
  $\,\star = \star_f \leq \star_a\,$ (Proposition~\ref{pr:4.1} (1)) 
  then,
  in
  particular, $\,N_{0}\,$ is also a quasi--$\star$--prime ideal of
  $\,D\,$.

  (2) Each quasi--$\star_a$--ideal of $\,D\,,$\ like $\,P\,,$\ is
  contained in a
  quasi--$\star_a$--maximal ideal $\,Q\,$ of $\,D\,$  (Lemma~\ref{lm:2}
  (1)).  In particular, we have
  $\,
   P \subseteq Q = Q^{\star_a} \cap D = Q{\textstyle\rm Kr}(D,\star) 
  \cap
   K \cap D = Q{\textstyle\rm Kr}(D,\star) \cap D\,.
  $
  \   Since $\,Q^{\star_a} \cap D \not = D\,,$ \ then necessa\-rily
  $\,Q{\textstyle\rm Kr}(D,\star) \neq {\textstyle \rm Kr}(D,\star)\,.$ 
  \
  Therefore there exists a maximal ideal of $\,{\textstyle\rm
  Kr}(D,\star)\,$ containing $\,Q{\textstyle\rm Kr}(D,\star)\,$ or,
  equi\-va\-le\-ntly, a valuation overring $\,(W,N)\,$ of
  $\,{\textstyle\rm Kr}(D,\star)\,$ with center in $\,{\textstyle\rm
  Kr}(D,\star)\,$ containing $\,Q{\textstyle\rm Kr}(D,\star)\,.$ \ By 
  (1,
  b), $N \cap D$ is a quasi--$\star_a$--prime ideal of $\,D\,.$ \ Since 
  it
  contains $\,Q\,,$\ by the maximality of $\,Q\,,$\ we deduce that $\,Q 
  =
  N \cap D\,$.

  (3) Since $\,{\textstyle\rm Kr}(D,\star)\,$ is a B\'ezout domain
  (Proposition~\ref{pr:4.1} (3)) then the maximal spectrum of
   $\,{\textstyle\rm Kr}(D,\star)\,$ is described by the centers in
  $\,{\textstyle\rm Kr}(D,\star)\,$ of the minimal valuation overrings 
  of
  $\,{\textstyle\rm Kr}(D,\star)\,.$ \ The conclusion follows 
  immediately
  from (2).

  (4) Recall that if $\,V\,$ is a $\,\star$--valuation overring of
  $\,D\,,$ the map
  $\,V \mapsto V(X)\,$ (where $\,V(X)\,$ is the trivial extension of
  $\,V\,$ into $\,K(X)\,$ \cite[page 218]{G}) defines an order 
  preserving
  bijection
  between the set of all the $\,\star$--valuation overrings of $\,D\,$ 
  and
  the set of all the valuation overrings of $\,{\textstyle\rm
  Kr}(D,\star)\,$ \cite[Theorem 3.5]{FL2}.  Therefore, by (3), if $\,Q 
  \in
  {\calM}(\star_a)\,,$\ we can find a (minimal) valuation overring
  $\,(W,N)\,$ of $\,{\textstyle\rm Kr}(D,\star)\,$ such that $\,N \cap 
  D =
  Q\,.$ \ Then, we can consider $\,V := W \cap K\,$ and $\,M := N \cap 
  K =
  N \cap V\,.$ \ By the previous remark, $\,V\,$ is a 
  $\,\star$--valuation
  overring of $\,D\,$ and its maximal ideal $\,M\,$ is such that $\,M 
  \cap
  D = Q\,.$ \ Hence, $\,(V,M)\,$ dominates $\,D_Q\,$. \hfill $\Box$

  %COROLLARY 4.4

  \begin{coor} \label{cor:4.4} \sl Let $\,\star\,$ be a nontrivial
  semistar
  operation of an integral domain $\,D\,$.  Assume that $\,\star\,$ is 
  an
  e.a.b. semistar operation
  of finite type with $\,D = D^\star\,.$ \
  Then each $\star$--maximal ideal of $\,D\,$ is the center in
  $\,D\,$ of a  minimal $\,\star$--valuation overring of $\,D\,.$

  \end{coor} \rm

  \vskip -3pt \noindent {\bf Proof}.  In the present situation, we have
  $\,{\calM}(\star_a) = \mbox{\rm Max}^\star(D)\,$.\, The conclusion
  follows easily from Theorem~\ref{thm:4.3} (3) and (4).  \hfill $\Box$

  %\vspace{.2in}

  %COROLLARY 4.5

  \begin{coor} \label{cor:4.5} \sl  Let $\,\star\,$ be a nontrivial
  semistar
  operation of an integral domain $\,D\,.$\ Assume that $\,\star =
  \star_f\,$. \ Then:

  %\begin{enumerate}
  \bara
  \rm \bf \item \sl $\tilde{\star} \leq \widetilde{(\star_a)} =
  (\star_a)_{sp}
  \leq \star_a\;\;$
   and $\;\;\tilde{\star} \leq (\tilde{\star})_a \leq \star_a \, $.

   \rm \bf \item \sl
  $ {\textstyle\rm Na}(D,\star) =
  {\textstyle\rm Na}(D,\tilde{\star}) \subseteq {\textstyle\rm
  Na}(D,\widetilde{(\star_a)}) = {\textstyle\rm Na}(D,\star_a) \subseteq
   {\textstyle\rm Kr}(D,\star_a) 
   = {\textstyle\rm Kr}(D,\star)\,.$

   \rm \bf \item \sl
   ${\textstyle\rm Na}(D,\star) = {\textstyle\rm
   Na}(D,\tilde{\star}) \subseteq {\textstyle\rm 
  Na}(D,(\tilde{\star})_a)
   \subseteq {\textstyle\rm Kr}(D,(\tilde{\star})_a) = {\textstyle\rm
  Kr}(D,\tilde{\star}) \subseteq {\textstyle\rm Kr}(D,\star) \,.$

  \rm \bf \item \sl For each $\,E \in \boldsymbol{\overline{F}}(D)\,,$

  %\begin{enumerate}
  \balf

  \rm \bf \item \sl $E^{\widetilde{ (\star_a}) } =
  E{\textstyle\rm Na}(D,\star_a) \cap K \; (\ \supseteq E{\textstyle\rm
  Na}(D,\star) \cap K =E^{\tilde{\star}}\ )\,;$

  \rm \bf \item \sl $E^{(\tilde{\star})_a} = E{\textstyle\rm
  Kr}(D,\tilde{\star})
  \cap K \; (\ \subseteq E{\textstyle\rm Kr}(D,\star) \cap K =
  E^{\star_a}\ )\,.$

  %\end{enumerate}
  \ealf
  %\end{enumerate}
  \eara
  \end{coor}

  \vskip -3pt \noindent {\bf Proof}.  (1) follows trivially from the 
  fact that if
  $\,\star_1\,$ and $\,\star_2\,$
  are two semistar operations of finite type, then
  $
  \star_1 \leq \star_2 \;$ implies that $\; (\star_1)_{sp} = 
  \widetilde{\
  \star_1\ } \leq
  \widetilde{\ \star_2\ } = (\star_2)_{sp}\,,
  $
  \ and from Corollary~\ref{cor:9} (1) and
  Proposition~\ref{pr:4.1} (1).

  (2) and (3) are consequences of Corollary~\ref{cor:3.5} (3) and of
  Proposition~\ref{pr:4.1} (2).

  (4) follows from Proposition~\ref{pr:3.3} (3) and 
  Proposition~\ref{pr:4.1} (5).  \hfill $\Box$

  \vspace{4pt}
% S E C T I O N
\section{The semistar operations \ $\tilde{\star}$ \ and \ $\star_a$}

  In this section we consider more closely the two operations which are
  naturally associated with the semistar Nagata rings and the semistar
  Kronecker
  function rings:

  % \begin{itemize}

  % \;\;\; $\bullet$ \, \;
  \centerline{ $\tilde{\star}\;$\ associated with \,Na$(D,\star) $ \;\; 
  \;  and \;\; \;
  $\star_a\;$ associated with \,Kr$(D,\star)\,.$ }
  % \;\;\; $\bullet$ \, \;
  \smallskip

  % \end{itemize}

  \noindent An elementary first question to ask is whether the two 
  semistar
  operations are actually the same - or usually the same - or rarely the
  same.  Theorem~\ref{tm:3.8} indicates that for a semistar operation
  $\,\star\,$ on a domain $\,D\,$, the $\,\tilde{\star}$--valuation
  overrings of $\,D\,$ are all the valuation overrings of the
  localizations of $\,D\,$ at the primes in $\,{\calM}(\star_f)\,$.\ On
  the other hand, \cite[Proposition 3.3 and Theorem 3.5]{FL2} indicates
  that the $\,\star_a$--valuation overrings (or, equivalently, the
  $\,\star$--valuation overrings) of $\,D\,$ correspond exactly to the
  valuation overrings of the Kronecker function ring \,Kr$(D,\star)\,$.\
  It is easy to imagine that these two collections of valuation domains
  can frequently be different.  We consider several different examples.

  %EXAMPLE 5.1

  \begin{exxe} \label{ex:5.1} \sl A (semi)star operation $\,\star\,$ of
  an integral domain $\,D\,$ such that $\,\tilde{\star} \neq
  \star_{a}\,,$\ but the $\,\tilde{\star}$--valuation overrings  
  coincide
  with the $\,\star_a$--valuation overrings (and so \
  \rm Kr$(D,\tilde{\star})=$ Kr$(D,\star_{a})=$ Kr$(D,\star),$\
  \cite[Corollary
  3.8]{FL1} \sl and \rm \cite[Theorem  3.5]{FL2}\sl).

  \rm  Let $\,L\,$ be a field and let $\,D\,$ be the localization
  $\,L[X,Y]_M\,$ of the polynomial ring  $\,L[X,Y]\,$ at the maximal 
  ideal
  $\,M :=(X,Y)\,.$ \ Let $\,\star := d\,$ be the identity (semi)star
  operation on
  $\,D\,$ (defined by $\,E^d := E\,,$ for all
  $\,E \in \overline{F}(D)\,).$\ Clearly, $\,\star_{f}=d=\star\,$ and
  every prime ideal of $\,D\,$ is
  a $\,\star_f$--prime. It follows that $\,\tilde{\star} = \star =
  d\,$ (Corollary~\ref{cor:9} (2, b)).\   On the other hand, note that 
  in general:

  \noindent \bf Claim. \sl For each integral domain $\,D\,,$\ the e.a.b.
  semistar operation $\,d_{a}\,$ as\-so\-cia\-ted to the identity
  (semi)star
  operation $\,d\,$ of $\,D\,$ coincides with the (semi)star operation 
  $\,b\,$
  defined, for each $\,E \in \overline{F}(D) ,\,$ by $ E^{b}:= \cap\{EV 
  \,|
  \; V \,\mbox{is a valuation overring of}\ D\}\,.  $ \rm

  \noindent The claim follows from Proposition~\ref{pr:4.1} (5) and
  \cite[Theorem  3.5]{FL2},
   since:
  $$\begin{array}{rl}
  E^{d_{a}} \hskip -9pt &= E\mbox{\rm Kr}(D, d) \cap K = \cap\{EV(X) \,|
   \,\, V \ \mbox{is a valuation overring of}\ D\} \cap K =\\
   \hskip -9pt &= \cap\{EV \, | \,\, V \,\mbox{is a valuation overring
  of}\ D\} = E^{b} \,.
   \end{array}
   $$

   Note that $\, d \neq b = d_{a} \,$ in $\,D\,$ (otherwise
  $\,D\,$
  would be a Pr\"ufer domain by \cite[Theorem 24.7]{G}), hence  
  $\,\star =
  d\,$ is not e.a.b.
  \cite[Proposition 4.5 (5)]{FL1} and so $\,\tilde{\star}\ (= \star = d)
  \neq
  \star_a\ (= d_{a} =b)\,.$\ Moreover, every valuation overring of 
  $\,D\,$ is
  (obviously) a $\,\tilde{\star}$--valuation overring and also (by the
  claim) every valuation overring of $\,D\,$ is a $\,\star_a$--valuation
  overring of $\,D\,.$\ Therefore, $\,\tilde{\star}\,$ and $\,\star_a\,$
  are different, but have the same collection of ``associated'' 
  valuation
  overrings.  \noindent Finally, observe that, for the particular
  $\,\star\,$ we are considering here, we have (using a ``new''
  indeterminate $\,Z\,$): $$ \begin{array}{rlr} \mbox{\rm Na}(D, \star)
  \hskip -5pt &=\mbox{\rm Na}(D, \tilde{\star}) =\mbox{\rm Na}(D,d) 
  =D(Z)
  \subsetneq & \\
  \hskip -5pt &\subsetneq \mbox{\rm Kr}(D, \star) =\mbox{\rm Kr}(D,
  \tilde{\star})= \mbox{\rm Kr}(D, \star_{a}) =\mbox{\rm Kr}(D, d) = & 
  \\
  \hskip -5pt &= \cap \{V(Z) \; | \; \, V \, \mbox{is an overring of} \,
  D\}\,.  & \quad \quad \quad \quad \Box
  \end{array}
  $$
  %  It follows
  %that $\,\star_a\,$ is actually the classical $\,b\,$ operation.
  \end{exxe}

  We noted in the preceding example that although $\,\tilde{\star}\,$ 
  and
  $\,\star_a\,$ were different, nevertheless, the collection of the
  $\,\tilde{\star}$--valuation overrings  of $\,D\,$ coincides with the
  collection of the $\,\star_a$--valuation overrings of $\,D\,.$\ The 
  next
  example displays wider differences between the two operations.

  %EXAMPLE 5.2

  \begin{exxe} \label{ex:5.2} \sl A (semi)star operation $\,\star\,$ of
  an integral domain $\,D\,$ such that $\,\tilde{\star} \neq
  \star_{a}\,,$ \ the $\,\star_a$--valuation overrings form a proper
  subset of the set of $\,\tilde{\star}$--valuation overrings, but
   $\,\tilde{\star} =\widetilde{(\star_a)}\,.$

   \rm Let $\,L\,$ and $\,D\,$ be as in Example~\ref{ex:5.1}.
   Let
  $\,N := MD\,$ denote the maximal ideal of $\,D\,$. \
   For each irreducible
  polynomial $\,f \in M\,,$\  let $\,W_f :=
  L[X,Y]_{(f)} = D_{(f)}\,.$ \ Then $\,W_f\,$ is a DVR overring of
  $\,D\,.$ \
  Let $\,V_{\mbox {\tiny \it X}}\,$ be the two dimensional valuation
  o\-ver\-ring of
  $\,D\,$ with maximal ideal generated by $\,Y\,$ and with
  $\,W_{X}=\,L[X,Y]_{(X)}\,$ as a one dimensional (valuation) overring,
  i.e.

  \centerline{$ V_{\mbox {\tiny \it X}} := L[Y]_{(Y)} + XL[X,Y]_{(X)} \
  (\subsetneq W_{X})\,.
  $}

   \noindent We consider the following family of valuation overrings of
  $\,D\,$:

  \centerline{$ {\calW}:= \{ W_f \;|\;\,f \in M\,, \;\, f \neq X\,,\; 
  f\ \;
  \mbox{irreducible in} \; L[X,Y]\} \, \cup \, \{V_{\mbox {\tiny \it
  X}}\}\,
  $}

  \noindent  and we define a semistar operation $\,\star :=
  \star_{\calW}\,$ of $\,D\,$ by setting $\,E^\star := \cap \{ EW \; |$
   $\;\, W \in {\calW} \}\,$ for all $\,E \in \overline{F}(D)\,.$ \ It 
  is
  well known that $\,\star_{\calW}\,$ is an e.a.b. (in fact, a.b.)
  semistar operation of  $\,D\,$  and the
  Kronecker function ring associated with $\,\star\,$ (in the ``new''
  variable
  $Z$) is then \ Kr$(D,\star) = \cap \{ W(Z) \;|\;\, W \in {\calW}
  \}\,$ \cite[Corollary 3.8]{FL2}.  We claim that:

  \noindent \bf Claim. \sl The maximal ideals of $\,\mbox{\rm
  Kr}(D,\star)\,$ are exactly the centers
   of the maximal ideals of the valuation domains  $\,W(Z)\,,$\  when 
  $\,W
  \in{\calW}\,$. \rm

   To prove the claim note that, from the fact that  \
  Kr$(D,\star)\,$ is a Pr\"ufer (in fact,
  B\'ezout) domain (Proposition~\ref{pr:4.1} (3)) and from \cite[Theorem
  3.5]{FL1},
  %each maximal ideal of\ Kr$(D,\star)\,$ must
  %contain one of the polynomials in $\,{\calW}\,$ (since
   there exists a canonical
  bijection between the  maximal ideals of\ Kr$(D,\star)\,$   and the
   valuation overrings of\ Kr$(D,\star)\,$ of the type $\,V(Z)\,,$\
   where $\,V\,$ is a minimal $\,\star$--valuation overring
  of $\,D\,$ (cf. also \cite{Dobbs-Fontana}). Moreover observe that, by
  definition, each $\, W \in {\calW}\,$ is a
   $\,\star$--valuation overring and that  the intersection $\,\cap \{
  W(Z)\; |\;\, W \in {\calW} \}\,$
   is irredundant, i.e., if any one of the valuation domains $\,W\,$ was
  omitted,
    the intersection would be different (in fact, it is easy to see that
    the first intersection in the following formula
    $$
    \begin{array}{rl}
    D \hskip -5pt &= \cap \{ W_f  \;|\;\,f \in M\,, \;\, f\ \;
  \mbox{irreducible
  in} \; L[X,Y]\}) =\\
  \hskip -5pt&= \cap \{ W_f  \;|\;\,f \in M\,, \;\, f \neq X\,,\; f\ \;
  \mbox{irreducible
  in} \; L[X,Y]\} \,
  \cap \, \{V_{\mbox {\tiny \it X}}\} =\\
  \hskip -5pt &=: \cap \{ W \; |\;\, W \in {\calW} \}
  \end{array}
  $$
  \noindent is irredundant, because
    $\,D\,$ is a Krull domain, so it is the same for the last
    intersection; this property implies easily the
    irredudancy of the $\,\cap \{ W(Z)\; |\;\, W \in {\calW} \}\,$).

    Note that the family of
    valuation overrings  $\,{\calW}(Z) := \{ W(Z)\; |$ $\;\, W \in 
  {\calW} \}\,$ of the  Pr\"ufer
    domain \ Kr$(D,\star) \,$ has finite character (in the sense
    that each nonzero element of  \ Kr$(D,\star) \,$ is a nonunit in
    at most finitely  many valuation overrings of $\,{\calW}(Z)\,$).

    \noindent As a matter of fact, in this case $\,\star = \star_{a}\,$
  hence $
    (a_{0}, a_{1}, \ldots, a_{n})^\star = f\mbox{Kr}(D, \star) \cap K\,,
  \; \,
    \mbox{for each}\; 0 \neq f := \sum_{k=0}^{n}a_{k}Z^k\in D[Z]\,,
    $ \ (cf. \cite[Theorem 3.11 (1), (2) and Theorem 5.1 (2)]{FL2})
  and
    $\,(a_{0}, a_{1}, \ldots, a_{n})^\star \in \boldsymbol{f}(D)\,,$ \
    because $\,D\,$ is a Noetherian ring; moreover,  each nonzero 
  finitely
  generated ideal of $\,D\,$ is contained in at most finitely many
  height
    1 prime ideals of $\,D\,,$ \ because $\,D\,$ is a Krull domain.

    The claim then follows immediately from \cite[Corollary
  1.11]{GH}. \ In other words, each maximal ideal\ $H$\ of \ $\mbox{\rm
  Kr}(D,\star)$ \ contracts onto a prime ideal of $\,D\,$ and thus
  contains a poly\-no\-mial $\,f \in M\,, \ f$ \ irreducible in
  the polynomial ring  $\, L[X,Y]\,;$ \ henceforth, if $\,f \neq X\,,$ \
  then
  $\, H = fW_{f}(Z) \cap \mbox{\rm Kr}(D,\star)\,,$ \ if $\,f = X\,,$ \
  then $\, H = YV_{\mbox {\tiny \it X}}(Z) \cap \mbox{\rm 
  Kr}(D,\star)\,.$

   The import of this claim is that the collection $ \{ W_f \,|\,\,f \in
   M\,,\, f \neq X\,,\, f \,$ irreducible in $ \; L[X,Y]\} \ \cup \
   \{V_{\mbox {\tiny \it X}}, \ W_{X} \}= {\calW} \ \cup \ \{W_{X} \}\ 
  , $
   constitutes the collection of all nontrivial $\,\star_a$--valuation 
  (or,
   equivalently, $\,\star$--va\-lua\-tion) overrings of $\,D\,.$

  On the other hand, note that $\,D\,$ is local and the maximal
  ideal of the valuation overring $\,V_{\mbox \small X} \in {\calW} \,$
   is centered on the maximal ideal $\,N\,$ of $\,D\,$.
    Moreover, $\,D\,$ is Noetherian, so $\,\star\,$ is a semistar
    operation of finite type on $\,D\,$. \ It follows that the maximal
    ideal $\,N\,$ of $\,D\,$
    belongs to $\,{\calM}(\star_f)\,.$ \ Since (obviously) $\,D_N = D\,$
  this leads
    to the conclusion that $\,{\calM}(\star_f) = \{N\} =
  \mbox{Max}(D)\,,$\ and
    so $\,\tilde{\star} =
  d\,$ (the identity (semi)star operation) of $\,D\,.$
   \ As noted in the previous example, this implies that every valuation
  overring of
   $\,D\,$ is a $\,\tilde{\star}\,$--valuation overring of $\,D\,$.
   \noindent Therefore, $\,\star_a\,$ and $\,\tilde{\star}\ (=d\ )\,$ 
  are
  not only
   different as semi\-star operations,
  but they are also associated with different sets of valuation
  overrings (e.g. for each  $\ f\in M\ ,\, f \neq X\ ,\  f \ $ 
  irreducible
  in  $\, L[X,Y]\,,$ \ the two dimensional valuation
  o\-ver\-ring $\,V_{\mbox {\tiny \it f}}\,$ of $\,D\,,$ \ having
  $\,W_{f}=\,L[X,Y]_{(f)}\,$ as a one dimensional (valuation) overring 
  and
  dominating $\,D\,,$ \ is a valuation overring of $\,D\,,$ \ but is 
  not a
  $\,\star_a$--valuation overring of $\,D\,$).

  In the present situation, observe that we have that $\, d =
  \tilde{\star} =\widetilde{(\star_a)}\,.$\ As a matter of fact, from
  Theorem~\ref{thm:4.3} (3) and from the fact that \ Kr$(D,\star)\,$ is 
  a
  Pr\"ufer domain, we have that each member of  $\,{\calM}(\star_a)\,$ 
  is the
  center in $\,D\,$ of a  minimal $\,\star$--valuation overring
  of $\,D\,,$\ thus $\,{\calM}(\star_a)= \{N\}\,,$\ hence
  $\,{\calM}(\star_a)= {\calM}(\star_{f})\,.$ \ Finally, we have:
  \vskip -12pt$$
  \begin{array}{rll} \mbox{\rm Na}(D, \star) \hskip -5pt &=\mbox{\rm 
  Na}(D,
  \tilde{\star}) = \mbox{\rm Na}(D, \widetilde{(\star_{a})})= \mbox{\rm
  Na}(D, \star_{a}) = \mbox{\rm Na}(D,d) =D(Z) \subsetneq & \\
      \hskip -5pt &\subsetneq \mbox{\rm Kr}(D, d) =
  \mbox{\rm Kr}(D, b) = \cap \{V(Z) \; | \; \, V \, \mbox{is an overring
  of} \, D\}\,\subsetneq & \\
        \hskip -5pt &\subsetneq \mbox{\rm Kr}(D, \star) = \mbox{\rm 
  Kr}(D,
  \star_{a}) = \cap \{W(Z) \; | \; \, W\, \in {\calW}\}\,.  & 
  \quad\quad \Box
  \end{array}
  $$
   %= \mbox{\rm Kr}(D, \tilde{\star}) = \mbox{\rm
  % Kr}(D,\widetilde{(\star_{a})})
  \end{exxe}

  %\vskip -12pt \hfill $\Box$

  \vskip 2pt

  In Example~\ref{ex:5.2} $\,\star_a \neq \tilde{\star}\,,$ \ however
  $\,\widetilde{(\star_a)} = \tilde{\star}\,.$ \ It seems plausible that
  something of this type holds in
  general.  The next example demonstrates that it does not and 
  illustrates
  why.

  %EXAMPLE 5.3

  \begin{exxe} \label{ex:5.3} \sl A (semi)star operation $\,\star\,$ of
  an integral domain $\,D\,$ such that $\,\tilde{\star} \neq
  \star_{a}\,,$ \ the $\,\star_a$--valuation overrings form a proper
  subset of the set of $\,\tilde{\star}$--valuation overrings and
  $\, \tilde{\star} \neq \widetilde{(\star_a)}\,.$

  \noindent \rm  Let $\,D\,$ and $\,N\,$ be as in the two previous
  examples.
  \ We construct a (semi)star
  operation $\,\star\,$ on $\,D\,$ as follows:

  \begin{enumerate}

  \rm \bf \item  \rm If $\,dD\,$ is any nonzero principal ideal of
  $\,D\,,$ \ then
  $\,(dD)^\star := dD\,$.

  \rm \bf \item  \rm  If $\,J \subseteq D\,$ is a nonzero ideal of 
  $\,D\,$
  which is not contained in any proper principal ideal of $\,D\,$, \ 
  then
  $\,J^\star := N\,$.

  \rm \bf \item  \rm If $\,J \subseteq D\,$ is a nonzero ideal of 
  $\,D\,$
  which is not principal,
   but is contained in a principal ideal, then we factor $\,J\,$ as 
  $\,J =
  fI\,,$ \ where
  $\,f\,$ is a GCD of a set of generators of $\,J\,$ and $\,I := (J
  :_{D} fD)\,$ is not contained
  in any proper principal ideal of $\,D\,$\ by the choice of $\,f\,.$ \
  Then
  $\,J^\star := fN\,$.

  \rm \bf \item  \rm  If $\,J\,$ is a nonzero fractionary ideal of 
  $\,D\,$
  which is
  not contained in $\,D\,$,\ choose a nonzero element $\,d \in D\,$ such
  that $\,dJ \subseteq D\,$. \
  Then define
  $\,J^\star := (1/d)(dJ)^\star\,$.

  \rm \bf \item  \rm If $\,J \in \overline{\boldsymbol{F}}(D) \setminus
  \boldsymbol{F}(D)\,$ we define $\,J^\star := L(X,Y)\,$.

  \end{enumerate}

  Since $\,D\,$ is Noetherian, then $\,\star\,$ is of finite
  type.  Henceforth, it is clear  that $\,{\calM}(\star_f) = \{N\}\,.$  
  \
  Thus, as in the previous example, $\,\tilde{\star} = d\,$.

  \noindent However, since $\,D\,$ is integrally closed and Noetherian, 
  it
  is easy to see from the definition
  of $\,\star_a\,$ that $ D^{\star_a} =D \,, \;\, \,N^{\star_a} = D\, 
  \;\,
  \mbox{and} \;\, (fD)^{\star_a} =fD\,,
  $ \ for each $\,  f \in M\,$ and $\,f\,$ irreducible
  in $L[X,Y]\,.$ \ Hence,
  $\,{\calM}(\star_a) = \{ fD \;| \;\ f \in M\,$ and $\,f\,$ irreducible
  in $L[X,Y]$ $\}\,$ coincides with the set of all the height 1 primes
  of $\,D\,.$ \ Moreover, $\,d_a\,$ is the classical $\,b\,$ operation
  (Claim in Example~\ref{ex:5.2}) and thus $\, (\tilde{\star})_{a} =
  d_{a} =b\,.$

   On the
  other hand,
  $\,\star_a\,$ coincides with the $\,t\,$ (semi)star operation of
  $\,D\,.$ \ As a matter of fact, we observed already that
  $\,{\calM}(\star_a)\,$ coincides with the set of all the height 1 
  primes
  of $\,D\,,$ \ and this implies that $\,\widetilde{(\star_a)} = t\,$
  because $\,D\,$ is a Krull domain \cite[Proposition 44.13 or Theorem
  44.2]{G}.  Since for (semi)star operations of finite type we have 
  always
  the inequalities $\,\widetilde{(\star_a)} \leq \star_a \leq t\,$
  (Corollary~\ref{cor:4.5} (1) and \cite[Theorem 34.1 (4)]{G}), we 
  deduce
  immediately that $\,\widetilde{(\star_a)} = \star_a = t\,$.

  Observe that $\,{\calM}(b) = \mbox{Max}(D)\,,$ since every
  valuation overring is a $\,b$--valuation overring and every prime 
  ideal
  of $\,D\,$ is a $\,b$--prime of $\,D\,.$ \ We conclude, for the
  particular $\,\star\,$ we are considering here, that:
   $$
  b= (\tilde{\star})_{a} \neq \widetilde{(\star_a)} =t\,,\;\;\; 
  \mbox{and}
  \;\;\;
  \, d = \tilde{b} = \widetilde{(\tilde{\star})_{a}} \neq
  \widetilde{\widetilde{(\star_a)}} = \widetilde{(\star_a)} =t\,.
  $$
   \noindent So it is hopeless to try to attain an equality by applying
   $\,\widetilde{(\mbox{-})}\,$ and $\,(\mbox{-})_a\,$ in different
  orders.

   \noindent Finally, observe that, for the particular $\,\star\,$ we 
  are
  considering here, we have (using a ``new'' indeterminate $\,Z\,$):
  $$ \begin{array}{rl}
       \mbox{\rm Na}(D, \star) \hskip -5pt &=\mbox{\rm
  Na}(D, \tilde{\star}) = \mbox{\rm Na}(D,d) =D(Z) \subsetneq \\
       \hskip -5pt &\subsetneq \mbox{\rm Na}(D, 
  \widetilde{(\star_{a})})=
    \mbox{\rm Na}(D, \star_{a}) = \mbox{\rm Na}(D, t) = \mbox{\rm Na}(D,
  v) =\\
     \hskip -5pt &=\cap \{W_{f}(Z) \; | \; \, f\in M\,, \;\,  f \; 
  \mbox{
  irreducible
  in} \, L[X,Y]\,\} =\\
       \hskip -5pt &= \mbox{\rm Kr}(D, \star)  =\mbox{\rm Kr}(D,
  \widetilde{(\star_{a})})
    =\mbox{\rm Kr}(D, \star_{a}) = \\
    \hskip -5pt &= \mbox{\rm Kr}(D, t) = \mbox{\rm Kr}(D, v) =
     \cap \{W(Z) \; | \; \, W\, \in {\calW}\}\,.\\
   
     \mbox{\rm Na}(D, \star) \hskip -5pt &=\mbox{\rm Na}(D,
    \tilde{\star}) = \mbox{\rm Na}(D,d) =D(Z) \subsetneq \\
    \hskip -5pt &\subsetneq \mbox{\rm Kr}(D, d) = \mbox{\rm Kr}(D,
  \tilde{b}) = \mbox{\rm Kr}(D, b) =\\
    \hskip -5pt &= \cap \{V(Z) \; | \; \, V \, \mbox{is an overring
    of} \, D\}\,\subsetneq \\
        \hskip -5pt &\subsetneq \mbox{\rm Kr}(D, \star) = \mbox{\rm 
  Kr}(D,

        \star_{a}) = \cap \{W(Z) \; | \; \, W\, \in {\calW}\}\,.
  \end{array}
  $$

  \end{exxe}

  \vskip -0.15in \hfill $\Box$

  \vskip 3pt

  It is possible to make a positive statement about the relationship
  between
  $\,\widetilde{(\mbox{-})}\,$ and $\,(\mbox{-})_a\,$ under conditions
  made clear in the
  preceding example.

  %PROPOSITION 5.4
  \begin{prro} \label{pr:5.4}  \sl Let $\,\star\,$ be a semistar
  operation of an integral domain $\,D\,.$ \
  Then, the following conditions are equivalent

  \brom

  \rm\bf \item \sl $\,\tilde{\star} = \widetilde{(\star_a)}\,;$

  \rm\bf \item \sl $\,{\calM}(\star_f) = {\calM}(\star_a)\,;$

  \rm\bf \item \sl \rm Na$(D, \star) = $ Na$(D, \star_{a})\,.$ \erom
  \end{prro} \rm

  \vskip -3pt \noindent {\bf Proof}.  (ii) $\Rightarrow$ (i).  If 
  ${\calM}(\star_f)
  = {\calM}(\star_a)\,$, then $\,\tilde{\star} = 
  \widetilde{(\star_a)}\,$
  follows immediately from the definition of the
  $\,\widetilde{(\mbox{-})}\,$
  operator and the fact that $\,\star_a\,$ is always (by definition) of
  finite type.

   (i) $\Rightarrow$ (ii).  If $\,{\calM}(\star_f) \not =
   {\calM}(\star_a)\,$, \ then $\,\tilde{\star} \not =
   \widetilde{(\star_a)}\,$ by Proposition~\ref{pr:3.3} (2) and (3), 
  again
   taking into account the fact that $\,(\star_a)_f = \star_a\,$.

   (iii) $\Rightarrow$ (i) and (i) $\Rightarrow$ (iii) follow from 
   Proposition~\ref{pr:3.1} (5) and Corollary~\ref{cor:3.5}.  \hfill
   $\Box$

  %\vspace{.2in}

  %REMARK 5.5
  \begin{reem} \rm Let $\,\star\,$ be a semistar
  operation of an integral domain $\,D\,.$ \ For semistar Kronecker
  function rings, we can easily state a result ``analogous'' to
  Proposition~\ref{pr:5.4} and concerning $\,\star\,$ and
  $\,\tilde{\star}\,.$ \
  More precisely, from
  Proposition~\ref{pr:4.1} (5), \cite[Theorem 3.1]{FL1} and
  \cite[Corollary 3.8, Theorem 5.1 (3)]{FL2}, we have that \sl the
  following conditions are equivalent:

  \brom

  \rm\bf \item \sl $\,\star_a = (\tilde{\star})_{a}\,$;

  \rm\bf \item \sl the set of $\,\tilde{\star}$--valuation
  overrings $\,D\,$ coincides with the set of  $\,\star$--valuation
  overrings
  of $\,D\,$;

  \rm\bf \item \sl $ \mbox{\rm Kr}(D, \tilde{\star}) = \mbox{\rm Kr}(D,
  \star) \,.$

  \noindent \sl\hskip -0.35in Moreover, each of the previous conditions 
  implies

  \rm\bf \item \sl $\,{\calM}(\star_a) = 
  {\calM}((\tilde{\star})_{a})\,.$

  \erom\rm
  \noindent On the other hand, observe that, from Example~\ref{ex:5.2}
  (for
  $\,\star = \star_{{\calW}}\,),$ \ we have that (iv)
  $\,\not\Rightarrow\,$ (iii), since $\,{\calM}(\star_a) = 
  {\calM}(\star_f) =
  \{N\}= {\calM}(b) ={\calM}(d_a) ={\calM}((\tilde{\star})_a)\,$ and 
  $\, \mbox{\rm Kr}(D, \tilde{\star}) =
  \mbox{\rm Kr}(D, b) \subsetneq  \mbox{\rm Kr}(D, \star)\,.$
  \end{reem}

  \vspace{2pt}

  This line of thinking motivates our final result,
   tying our investigation of diffe\-rently constructed semistar 
  operations
  back to the topic of
  Nagata rings.

  %PROPOSITION 5.5
  \begin{prro} \label{pr:5.5} \sl Suppose $\,\star_1\,$ and 
  $\,\star_2\,$
  are semistar operations on a domain $\,D\,$.  \ Then,
  $\,\mbox{\rm Na}(D,\star_1) = \mbox{\rm
  Na}(D,\star_2)\,$ if and only if
   $\,{\calM}((\star_1)_f) = {\calM}((\star_2)_f)\,.$
  \end{prro} \rm

  \vskip -3pt \noindent {\bf Proof}. First, suppose $\,\mbox{\rm 
  Na}(D,\star_1) =
  \mbox{\rm Na}(D,\star_2)\,$. \
  Then $\,{\calM}((\star_1)_f) = {\calM}((\star_2)_f)\,$ follows from
  Proposition~\ref{pr:3.1} (5).

  Now suppose that
  $\,{\calM}((\star_1)_f) = {\calM}((\star_2)_f)\,$. \
  Then $\,\mbox{\rm Na}(D,\star_1) = \mbox{\rm Na}(D,\star_2)\,$ follows
  from Proposition~\ref{pr:3.1} (4). \hfill  $\Box$

  \vspace{6pt}

  \noindent \bf Acknowledgements.  \rm During the preparation of this 
  work
  the first named author was partially supported by a research grant 
  MIUR
  2001/2002 (Cofin 2000 - MM 01192794).  Fontana thanks also the
  Mathematical Department of the Ohio State University for the 
  hospitality
  accorded during his visit in July 2001.

   \end{document}